\documentclass[10pt,oneside,a4paper]{amsart}
\usepackage{amsmath,amsfonts,amsthm, amssymb}
\usepackage{array}
\usepackage[all,textures]{xy}

\theoremstyle{plain}
\newtheorem{theorem}{Theorem}[section]
\newtheorem{proposition}[theorem]{Proposition}

\newtheorem{corollary}[theorem]{Corollary}

\theoremstyle{definition}
\newtheorem{definition}[theorem]{Definition}

\theoremstyle{remark}

\newtheorem{remark}[theorem]{Remark}

\newcommand{\Kern}{\mathrm{Ker}}

\renewcommand{\lim}{\mathrm{lim}}
\newcommand{\colim}{\mathrm{colim}}

\newcommand{\Tor}{\mathrm{Tor}}
\newcommand{\Ext}{\mathrm{Ext}}
\newcommand{\rExt}{\mathcal{E}\mathit{xt}}
\newcommand{\Hom}{\mathrm{Hom}}
\newcommand{\rHom}{\mathcal{H}\mathit{om}}
\newcommand{\rMod}{\mathcal{M}\mathit{od}}

\newcommand{\op}{^{\mathrm{op}}}
\newcommand{\Ob}{\mathrm{Ob}}

\newcommand{\AAA}{\mathfrak{a}}
\newcommand{\BBB}{\mathfrak{b}}

\newcommand{\OOO}{\mathfrak{o}}

\newcommand{\UUU}{\mathfrak{u}}
\newcommand{\VVV}{\mathfrak{v}}

\newcommand{\Alg}{\ensuremath{\mathsf{Alg}} }

\newcommand{\Ab}{\ensuremath{\mathsf{Ab}} }
\newcommand{\Mod}{\ensuremath{\mathsf{Mod}} }
\newcommand{\Fib}{\ensuremath{\mathsf{Fib}} }

\newcommand{\Sh}{\ensuremath{\mathsf{Sh}} }
\newcommand{\mmod}{\ensuremath{\mathsf{mod}} }

\newcommand{\Qch}{\ensuremath{\mathsf{Qch}} }
\newcommand{\Cat}{\ensuremath{\mathsf{Cat}} }

\newcommand{\Des}{\ensuremath{\mathrm{Des}}}

\newcommand{\lra}{\longrightarrow}

\newcommand{\aaa}{\ensuremath{\mathcal{A}}}
\newcommand{\bbb}{\ensuremath{\mathcal{B}}}
\newcommand{\ccc}{\ensuremath{\mathcal{C}}}
\newcommand{\ddd}{\ensuremath{\mathcal{D}}}

\newcommand{\fff}{\ensuremath{\mathcal{F}}}
\newcommand{\GGG}{\ensuremath{\mathcal{G}}}

\newcommand{\LLL}{\ensuremath{\mathcal{L}}}
\newcommand{\mmm}{\ensuremath{\mathcal{M}}}

\newcommand{\ooo}{\ensuremath{\mathcal{O}}}
\newcommand{\ppp}{\ensuremath{\mathcal{P}}}
\newcommand{\qqq}{\ensuremath{\mathcal{Q}}}

\newcommand{\sss}{\ensuremath{\mathcal{S}}}
\newcommand{\ttt}{\ensuremath{\mathcal{T}}}
\newcommand{\uuu}{\ensuremath{\mathcal{U}}}
\newcommand{\vvv}{\ensuremath{\mathcal{V}}}

\newcommand{\zzz}{\ensuremath{\mathcal{Z}}}
\CompileMatrices
\SelectTips{cm}{11}

\title{Algebroid prestacks and deformations of ringed spaces}
\author{Wendy Lowen$^{\ast}$}
\address{Departement DWIS\\ Vrije Universiteit Brussel\\ Pleinlaan
2\\1050 Brussel\\ Belgium/ Institut de math\'ematiques de Jussieu (UMR 7586)\\ 175, rue du Chevaleret\\ 75013 Paris }
\email[Wendy Lowen]{wlowen@vub.ac.be / lowen@math.jussieu.fr}
\thanks{$^{\ast}$Postdoctoral fellow FWO/CNRS}

\begin{document}
\begin{abstract}
 For a ringed space $(X,\ooo)$, we show that the deformations of the
  abelian category $\Mod(\ooo)$ of sheaves of $\ooo$-modules 
  \cite{lowenvandenbergh1} are obtained from algebroid prestacks, as
  introduced by Kontsevich.  In case $X$ is a quasi-compact separated
  scheme the same is true for $\Qch(\ooo)$, the category of
  quasi-coherent sheaves on $X$. It follows in particular that there is a
  deformation equivalence between $\Mod(\ooo)$ and $\Qch(\ooo)$.
\end{abstract}
\maketitle
\tableofcontents

\section{Introduction}

In \cite{kontsevich} Kontsevich proves that any Poisson bracket on a
$C^\infty$-manifold can be canonically quantized.  Similarly, in the
algebraic case one would like to quantize the structure sheaf $\ooo$ of a
smooth algebraic variety $X$. This can to a certain extent be done but
the gluing questions are more delicate and extra conditions are needed
(see \cite{yekutieli1}).

However in \cite{kontsevich1} Kontsevich takes a different approach:
corresponding to a Poisson bracket on $X$ he introduces a deformation 
of $\ooo$ in the category of 
\emph{algebroid prestacks} on $X$.  An algebroid prestack is the linear
analogue of a gerbe (see below for a precise definition).

The keypoint is that to an algebroid prestack one can associate
canonically an abelian category of coherent sheaves (in the noetherian
case).  One may think of this abelian category as a deformation of
$\operatorname{coh}(\ooo)$. Hence in this way the quantization of a
Poisson bracket in the algebraic case is achieved in complete
generality.

\medskip

In the current paper we show that Kontsevich's approach is very
natural and that indeed under weak hypotheses all deformations of the
abelian categories relevant to algebraic geometry are obtained from
 algebroids prestacks. We will do this in the framework of  the general
(infinitesimal) deformation theory of abelian categories which was
developed in \cite{lowenvandenbergh1}.  

\medskip

Here is the definition of an algebroid prestack \cite{polesello, kontsevich1}, which will
be a central notion in this paper. Recall first that a \emph{fibered category} \cite{moerdijk} is, roughly speaking, a presheaf of categories where the restriction functors commute only up to given isomorphisms, a \emph{prestack} is a fibered category which satisfies ``gluing for maps'' 
but not necessarily for objects, and a \emph{stack} \cite{moerdijk} is a prestack satisfying ``gluing for objects''. 
\begin{definition}
An \emph{algebroid prestack} $\aaa$ on a topological space $X$ is
a prestack of linear categories satisfying the following properties.
\begin{enumerate}
\item[(1)] Any point $x\in X$ has a neighborhood $U$ such that
  $\aaa(U)\neq \emptyset$.
\item[(2)] If $A,B \in \aaa(U)$ then every point $x\in U$ has a neigborhood
$V\subset U$ such that $A|_V\cong B|_V$. 
\end{enumerate}
\end{definition}
If $\ooo$ is a sheaf of rings on $X$ and if for $U\subset X$ open we
identify $\ooo(U)$ with a one-object category then $\ooo$ trivially
defines an algebroid prestack on $X$. Note that this is not a stack.
The associated stack is 
\begin{equation}
\label{verstacking1}
U\mapsto \{\text{locally free rank one $\ooo_U$-modules}\}
\end{equation}
Since we prefer to work with $\ooo$ rather than with \eqref{verstacking1} we will use algebroid prestacks rather than algebroid stacks.

\medskip 

Next we briefly sketch the theory developed in \cite{lowenvandenbergh1, lowenvandenbergh2}. We
consider deformations of abelian categories along a surjective map $S \lra R$ of
commutative coherent rings with nilpotent kernel $I$. An
$S$-deformation of an $R$-linear abelian category $\ccc$ is by
definition an $S$-linear abelian category $\ddd$ together with an
equivalence $\ddd_R\cong \ccc$ where $\ddd_R$ is the full subcategory of $R$-objects in $\ddd$, i.e. those objects
annihilated by $I$. In order to control this deformation
theory one has to restrict oneself to flat abelian categories. This is
a relatively technical notion but for an $R$-linear category with
enough injectives it simply means that the $\Hom$-sets between
injectives are $R$-flat.

One of the results of \cite{lowenvandenbergh1} is that for an $R$-algebra
$A$, there is an equivalence between
\begin{enumerate}
\item flat abelian deformations of $\Mod(A)$
\item flat algebra deformations of $A$
\end{enumerate}
Or in geometric terms: a deformation of an afine space is affine.

The key point to prove this result is that the finitely generated
projective generator $A$ of $\Mod(A)$ can be lifted uniquely (up to isomorphism) to any abelian
deformation.  This follows from the fact that when $I^2=0$ the liftings of $A$ are governed by an obstruction theory \cite{lowen2} involving $\Ext_A^{1,2}(A,I\otimes_R A)$, which are zero.  \medskip

On could hope for a similar equivalence between deformations of
$\Mod(\ooo)$ and $\ooo$ but this is only true if we consider
deformations of $\ooo$ in a more general category than ringed spaces.

To be able to state Theorem \ref{main1} below which describes the
deformations of $\Mod(\ooo)$ and which is one of our main results, we
note that for an algebroid prestack $\aaa$ on $X$ we may define an
associated abelian category $\Mod(\aaa)$. The objects of $\Mod(\aaa)$ are given by
the linear prestack maps $\aaa\lra
\mmm\mathit{od}(\underline{R})$, where
$\mmm\mathit{od}(\underline{R})$ denotes the stack of sheaves
of $R$-modules on $X$.  It is easy to see that this definition
gives the expected result for $\aaa=\ooo$. The following is contained in Theorem \ref{theorem}.
\begin{theorem}
\label{main1}
  Let $(X,\ooo)$ be a flat $R$-linear ringed space. Then every
flat $S$-deformation of $\Mod(\ooo)$ is of the form $\Mod(\aaa)$ where
  $\aaa$ is an $S$-linear algebroid prestack on $X$, which is a flat weak
deformation of $\ooo$ (in a sense to be made precise, see \S\ref{lindef}).
\end{theorem}

To get an idea how the algebroid prestack $\aaa$ may be
constructed, let $\ddd$ be a flat abelian $S$-deformation of
$\Mod(\ooo)$.  Using localization theory $\ddd$ may be transformed
into a deformation of the stack $\mmm\mathit{od}(\ooo)$ in the
category of stacks of abelian categories on $X$. Then we put
\begin{equation}
\label{dealgebroid}
\aaa(U)=\{\text{flat liftings of $\ooo_U$ to $\ddd(U)$}\}
\end{equation}
The fact that this prestack satisfies the conditions (1)(2) above
follows from the obstruction theory in \cite{lowen2}.

In Theorem \ref{theorem} we obtain a 1-1 correspondence between flat deformations of $\Mod(\ooo)$ and of $\ooo$ (in the sense explained in \S \ref{lindef}). The proof of this correspondence uses a ``liftable'' characterization of linear prestack maps $$\aaa \lra \ccc$$ inducing an equivalence of stacks $$\ccc \cong \rMod(\aaa).$$ This characterization is stated in the main Theorem \ref{karmodthm} of chapter 2. It is a ``local'' analogue of the standard characterization of linear maps $\AAA \lra \ccc$ inducing an equivalence of categories $\ccc \cong \Mod(\AAA)$ (Theorem \ref{baby}). Our proof of Theorem \ref{karmodthm} uses additive sheaf theory, and in particular the main result of \cite{lowen1}.

\medskip

If we now turn from an arbitrary ringed space to a ringed space $(X,\ooo)$ with an acyclic basis $\uuu$ (see \S \ref{paracyclic}), the situation of Theorem \ref{main1} becomes remarkably simpler. Indeed, starting from a flat abelian $S$-deformation $\ddd$ of $\Mod(\ooo)$ as above, the acyclicity condition implies, again using the obstruction theory of \cite{lowen2}, that for $U \in \uuu$ the object $\ooo_U$ can be lifted uniquely (up to isomorphism) to $\ddd(U)$. Hence $\aaa(U)$ defined in \eqref{dealgebroid} is itself an \emph{algebroid}, i.e. a nonempty linear category in which all objects are isomorphic. Consequently, $\aaa|_{\uuu}$ is equivalent to a \emph{twisted presheaf} on $\uuu$ (i.e. a presheaf where the restriction maps commute only up to given isomorphisms, see Definition \ref{twisted}). The following is contained in Theorem \ref{acyc}.

\begin{theorem}\label{main2}
Let $(X,\ooo)$ be a flat $R$-linear ringed space with acyclic basis $\uuu$. There is an equivalence between
\begin{enumerate}
\item flat deformations of the abelian category $\Mod(\ooo)$
\item flat deformations of $\ooo|_{\uuu}$ as a twisted presheaf
\item flat deformations of $\ooo|_{\uuu}$ as a fibered category (in the obvious, pointwise sense)
\end{enumerate}
\end{theorem}
Since deforming the fibered category $\ooo|_{\uuu}$ is readily seen to be equivalent to deforming the linear category $\OOO_{\uuu}$ with $\Ob(\OOO_{\uuu}) = \uuu$ and
$$\OOO_{\uuu}(V,U) = \begin{cases} \ooo(V) & \text{if}\, \, \, V \subset U \\ 0 & \text{else} \end{cases}$$ associated to $\ooo|_{\uuu}$ in \cite{lowenvandenbergh1}, Theorem \ref{main2} is actually a reformulation of 
\cite[Theorem 8.18]{lowenvandenbergh1}. An important advantage of this reformulation is that it allows us to make the connection with deformations of quasi-coherent sheaves.

\medskip

So, finally, let $(X, \ooo)$ be a quasi-compact separated scheme with a finite affine covering $\vvv$ which is closed under intersections and let $\uuu$ be the basis of all affine opens contained in some $V \in \vvv$. Suppose we are given a deformation $\aaa$ of the prestack $\ooo|_{\uuu}$. For every $U$, there is an equivalence between flat linear deformations of $\ooo(U)$ and flat abelian deformations of the module category $\Mod(\ooo(U))$. These deformations constitute a stack of deformed module categories $\Mod(\aaa(U))$ on $\uuu$, and we can ``glue'' them together to obtain a deformation $\Qch(\uuu, \aaa)$ of $\Qch(X)$. In fact it suffices to glue the categories $\Mod(\aaa(V))$ for $V \in \vvv$. The following final result is contained in Theorem \ref{quasi}. 
\begin{theorem}\label{finally}
Let $(X, \ooo)$ be a quasi-compact separated scheme with a finite affine covering $\vvv$ which is closed under intersections. There is an equivalence between
\begin{enumerate}
\item flat deformations of the abelian category $\Mod(\ooo)$
\item flat deformations of the abelian category $\Qch(\ooo)$
\item flat deformations of $\ooo|_{\vvv}$ as a twisted presheaf
\item flat deformations of the linear category $\OOO_{\vvv}$ associated to $\ooo|_{\vvv}$
\end{enumerate}
\end{theorem}
The analogue of Theorem \ref{finally} for Hochschild cohomology (which in particular yields Theorem \ref{finally} if $R$ is a field of characteristic zero) has been proved in \cite{lowenvandenbergh2}.

\section{Linear fibered categories and stacks of sheaves}

\subsection{Introduction}\label{intro1}

Throughout, $R$ will be a coherent, commutative ring. 

Sections \S \ref{fibgradcat}, \S \ref{parstacks}, \S \ref{parsheaves} of this chapter contain some preliminaries on fibered categories, prestacks, stacks and sheaves over them in the context of $R$-linear categories. Roughly speaking, an $R$-linear fibered category is a presheaf of $R$-linear categories where the restriction functors commute only up to given isomorphisms. Let $\aaa$ be an $R$-linear fibered category on a topological space $X$. Associated to $\aaa$ we have a stack $\rMod(\aaa)$ of sheaves on $\aaa$ and a canonical ``Yoneda-morphism'' $\aaa \lra \rMod(\aaa)$ of fibered categories. The aim of this chapter is to characterize this morphism intrinsically in a way that can be lifted under deformation. In the next chapter, we will use this characterization to prove that an abelian deformation of the stack $\rMod(\aaa)$ is again of the form $\rMod(\bbb)$ for a certain kind of deformation $\bbb$ of $\aaa$. In this chapter however, there is no reference to deformations (except to point out where a certain result will be used later on).

The characterization of $\aaa \lra \rMod(\aaa)$ (which is completed in \S \ref{parkarmod}) is twofold. First of all, we need to pinpoint some liftable properties of a stack of Grothendieck categories $\ccc$, which ensure that the restriction functors $j^{\ast}: \ccc(U) \lra \ccc(V)$ for $V \subset U$ are exact and come equipped with a fully faithful right adjoint $j_{\ast}$ and a fully faithful exact left adjoint $j_{!}$ (as is the case for $\rMod(\aaa)$). We give such conditions in terms of localizing subcategories in \S \ref{itoloc}. Next we need to characterize morphisms 
\begin{equation}
\aaa \lra \ccc
\end{equation} yielding an equivalence of stacks $\ccc \cong \rMod(\aaa)$. This will be done using additive topologies and sheaves. We develop the necessary preliminaries on this subject in \S \ref{parGFFF}. For every open $U$, we associate to the fibered category $\aaa|_U$ an additive category $\AAA_U$, which in the present setting comes with a natural morphism
\begin{equation}
\AAA_U \lra \ccc(U)
\end{equation}
The category $\AAA_U$ naturally inherits an additive Grothendieck topology $\ttt_U$ from the canonical topology on $X$, and we have an equivalence $$\Mod(\aaa|_U) \cong \Sh(\AAA_U, \ttt_U)$$
According to \cite{lowen2}, all we have to do to characterize (3) yielding $\ccc \cong \rMod(\aaa)$ is require that
\begin{itemize}
\item every morphism (4) satisfies the conditions (G), (F) and (FF) (see \S \ref{parGFFF})
\item the additive topology that $\AAA_U$ inherits from $\ccc(U)$ is precisely $\ttt_U$
\end{itemize}
Equivalently, we can require that 
\begin{itemize}
\item each $(\AAA_U, \ttt_U) \lra (\ccc(U), \ttt_{\mathrm{epi}})$ satisfies (G), (F) and (FF)
\item the objects of $\AAA_U$ become locally finitely presented and locally projective (Definition \ref{locep}) in $\ccc$.
\end{itemize}
This is precisely the statement of Theorem \ref{karmodthm}, which is thus a perfect analogue of the characterization of the Yoneda embedding $\AAA \lra \Mod(\AAA)$ of a linear category into its module category :  a  functor
\begin{equation}
\AAA \lra \ccc
\end{equation}
from $\AAA$ into a Grothendieck category $\ccc$ yields an equivalence $\ccc \cong \Mod(\AAA)$ prescisely when 
\begin{itemize}
\item (5) satisfies (G), (F) and (FF)
\item the additive topology on $\AAA$ induced by $\ttt_{\mathrm{epi}}$ on $\ccc$ is the trivial topology
\end{itemize}
or, equivalently, when
\begin{itemize}
\item (5) is fully faithful and the objects of $\AAA$ are generators in $\ccc$
\item the objects of $\AAA$ are finitely presented and projective in $\ccc$.
\end{itemize}

With the eye on application in the next chapter, we introduce some notions of flatness for fibered categories in \S \ref{parflat}, and we prove some preliminaries on algebroid fibered categories in \S \ref{algbfibcats}.

\subsection{Fibered graded categories}\label{fibgradcat}
For the classical theory of fibered categories we refer the reader to \cite{SGA1}. In this section we briefly present linear versions of some of the basic concepts. This involves the notion of a linear category graded over a base category. 

Let $\uuu$ be a base category and $R$ a commutative ring. A \emph{$\uuu$-graded $R$-linear category ($\uuu$-$R$-category) $\AAA$} consists of a (non-linear) category $\AAA$ and a functor $F: \AAA \lra \uuu$ such that:
\begin{itemize}
\item for every $A, A' \in \AAA$ and $f: F(A) \lra F(A')$ in $\uuu$, the fiber $\AAA_f(A,A') = \{ a \in \AAA \, | \, F(a) = f \}$ has an $R$-module stucture
\item for every $A, A', A'' \in \Ob(\AAA)$, $f: F(A) \lra F(A')$, $g: F(A') \lra F(A'')$, composition defines an $R$-module morphism $\AAA_g(A',A'') \otimes \AAA_f(A,A') \lra \AAA_{gf}(A,A'')$
\end{itemize}
A $\uuu$-$R$-category $\AAA$ has an \emph{associated $R$-linear category} $\AAA$ with the same object set and with $\AAA(A,A') = \coprod_{\uuu(F(A),F(A')}\AAA_f(A,A')$. For $U \in \uuu$, we denote by $\AAA(U)$ the fiber category of all objects $A$ with $F(A) = U$ and morphisms $a$ with $F(a) = 1_A$.

Note that if $\Ob(\uuu) = \Ob(\AAA) = \{ \ast \}$, then $\uuu$ is a semigroup and $\AAA$ is a $\uuu$-graded $R$-algebra in the classical sense.

There are obvious notions of $\uuu$-$R$-functors and $\uuu$-$R$-natural transformations making $\uuu$-$R$-categories into a 2-category $2\Cat(\uuu, R)$. A \emph{$\uuu$-$R$-functor} from $F: \AAA \lra \uuu$ to $G: \BBB \lra \uuu$ is a functor $K: \AAA \lra \BBB$ with $GK = F$ and defining $R$-module morphisms $K_f: \AAA_f(A,A') \lra \BBB_f(K(A),K(A'))$. A \emph{$\uuu$-$R$-natural transformation} $\eta: K \lra L$ between $K, L: \AAA \lra \BBB$ is an ordinary natural transformation with $\eta_A \in \BBB_1(K(A),L(A))$ (for $A \in \AAA$). The 2-category structure of $2\Cat(\uuu, R)$ yields a natural notion of \emph{equivalence of $\uuu$-$R$-categories}. Similar to the case of ordinary categories $K: \AAA \lra \BBB$ is an equivalence if and only if all the $K_f: \AAA_f(A,A') \lra \BBB_f(K(A),K(A'))$ are isomorphisms and all the $K_U: \AAA(U) \lra \BBB(U)$ are essentially surjective (hence equivalences of categories).

Let $\AAA$ be a $\uuu$-graded $R$-linear category.
A morphism $a \in \AAA_f(A,A')$ is called \emph{cartesian} if for every $B$ and $g: F(B) \lra F(A)$, the map $\AAA_g(B,A) \lra \AAA_{fg}(B,A')$ is an isomorphism. Note that for $U \in \uuu$, every section in $\AAA(U)$ is obviously cartesian. The $\uuu$-$R$-category $\AAA$ is called \emph{fibered} if for every $f: V \lra U$ in $\uuu$ and $A$ in $\AAA$ with $F(A) = U$, there is a $B$ with $F(B) = V$ and a cartesian morphism $a \in \AAA_f(B,A)$. A morphism (1-cell) of fibered $\uuu$-$R$-categories is an $\uuu$-$R$-functor preserving cartesian morphisms. Fibered $\uuu$-$R$-categories inherit the 2-category structure of $2\Cat(\uuu, R)$, yielding a 2-category $2\Fib(\uuu, R)$ together with a 2-functor $2\Fib(\uuu, R) \lra 2\Cat(\uuu, R)$.

Let $\AAA$ be a fibered $\uuu$-$R$-category. Put $\aaa(U) = \AAA(U)$.
Suppose we choose for every $A \in \AAA(U)$ and $f: V \lra U$ a cartesian morphism $f^{\ast}(A) \lra A$. Then clearly $\aaa(f) = f^{\ast}$ defines an $R$-linear functor $\aaa(U) \lra \aaa(V)$. This makes $\aaa$ into a \emph{pseudofunctor} from $\uuu^{\op}$ to the 2-category $2\Cat(R)$ of $R$-linear categories, i.e.
\begin{itemize}
\item an $R$-linear category $\aaa(U)$ for every $U \in \uuu$
\item an $R$-linear functor $i^{\ast}: \aaa(U) \lra \aaa(V)$ for every $i: V \lra U$ in $\uuu$
\item a natural isomorphism $\tau_{i,j}: (ij)^{\ast} \lra j^{\ast}i^{\ast}$ for each $j: W \lra V$, $i: V \lra U$
\end{itemize}
These data have to satisfy a ``cocycle condition'' for three composable morphisms, expressing that for an additional $k: Z \lra W$, the two canonical maps $(ijk)^{\ast} \lra k^{\ast}j^{\ast}i^{\ast}$ are identical.
Pseudofunctors from $\uuu^{\op}$ to $2\Cat(R)$ constitute a 2-category $\mathsf{Pseudo}(\uuu^{\op},2\Cat(R))$ (endowed with ``pseudo-natural transformations'' and ``modifications''). 

To a pseudofunctor $\aaa$ corresponds a fibered $\uuu$-$R$-category $\AAA$
with $$\Ob(\AAA) = \coprod_{\uuu}\Ob(\aaa(U))$$ and for $A_V \in \aaa(V)$, $A_W \in \aaa(W)$, $f: W \lra V$
$$\AAA_f(A_W,A_V) = \aaa(W)(A_W, f^{\ast}(A_V))$$
For $A_V \in \aaa(V)$ and $f: W \lra V$, a cartesian morphism in $\AAA$ is given by the identity morphism $f^{\ast}(A_V) \lra f^{\ast}(A_V)$ in $\aaa(W)$.
The correspondence can be made into a 2-equivalence $\mathsf{Pseudo}(\uuu^{\op},2\Cat(R)) \lra 2\Fib(\uuu, R)$.

In the sequel, we will use the term \emph{fibered category} interchangeably for a fibered $\uuu$-$R$-category or for a corresponding pseudofunctor (depending on the choice of cartesian morphisms).

If $\ooo: \uuu^{\op} \lra \Alg(R)$ is a presheaf of $R$-algebras (i.e. an honest functor), then $\ooo$ obviously defines a fibered category. The associated $R$-linear category $\OOO$ has $\Ob(\OOO) = \uuu$ and
$$\OOO(V,W) = \begin{cases} \ooo(V) & \text{if}\, \, \, V \subset W \\ 0 & \text{else} \end{cases}$$

\subsection{Fibered categories, prestacks and stacks on a topological space}\label{parstacks}

For an introduction to fibered categories and stacks on a topological space we refer the reader to \cite{moerdijk}. In this section we recall some of the basic concepts. 

Let $X$ be a topological space and let $\uuu \subset \mathrm{Open}(X)$ be a full subcategory of the category of open sets and inclusions.
In this case a $\uuu$-graded $R$-linear category is ``the same'' as an $R$-linear category $\AAA$ with $\AAA(V,U) = 0$ unless $V \subset U$.

As explained in the previous section, we will call a pseudofunctor $\fff$ from $\uuu$ to $2\Cat(R)$ an \emph{$R$-linear fibered category on $\uuu$}. For an inclusion $i: V \subset U$, we will often use the notation $i^{\ast}(F) = F|_V$.
We define the restriction $\fff|_U$ to $U \in \uuu$ to be the fibered category on $\uuu/U$ with $\fff|_U(V) = \fff(V)$.

Two objects $F, F' \in \fff(U)$ determine a presheaf $\Hom_{\fff}(F,F')$ on $\uuu/U$ with $$\Hom_{\fff}(F,F')(V) = \Hom_{\fff_V}(F|_V,F'|_V).$$ Suppose $\uuu$ is a basis of $X$. We denote the sheafication of the presheaf $\Hom_{\fff}(F,F')$ by $\rHom_{\fff}(F,F')$.
The fibered category $\fff$ is called a \emph{prestack} if the presheaves $\Hom_{\fff}(F,F')$ are sheaves.
If the categories $\fff(U)$ are abelian, we likewise define presheaves $\Ext_{\fff}(F,F')$ and their sheafications $\rExt_{\fff}(F,F')$.

For a fibered category $\fff$ and a covering $U_i \lra U$, there is an associated category of ``descent data'' $\mathrm{Des}(U_i,\fff)$ and a functor $\fff(U) \lra \mathrm{Des}(U_i,\fff)$. To define a \emph{descent datum}, we consider $\vvv = \{ V \in \uuu\, |\,\exists i\, V \subset U_i\}$. A descent datum consists of an object $F_V \in \fff(V)$ for every $V \in \vvv$, together with isomorphisms $F_V|_W \cong F_W$ for every inclusion $W \subset V$ in $\vvv$. These isomorphisms have to satisfy a compatibility condition for every two inclusions $W' \subset W \subset V$.  The fibered category $\fff$ is called a \emph{stack} if and only if, for every covering $U_i \lra U$, this functor is an equivalence of categories. 

If $\ooo$ is a sheaf of $R$-algebras on $\uuu$, then $\ooo$ is naturally a prestack but \emph{not} a stack. This is why we will continue to work with prestacks rather than stacks.

Let $\fff$ and $\GGG$ be $R$-linear fibered categories over $X$. A morphism of $R$-linear fibered categories (a ``pseudo-natural transformation'') $\phi: \fff \lra \GGG$ consists of the following data:
\begin{itemize}
\item an $R$-linear functor $\phi_U: \fff(U) \lra \GGG(U)$ for every $U \in \uuu$
\item a natural isomorphism $\alpha_i: \phi_Vi^{\ast} \lra i^{\ast}\phi_U$ for every $i: V \lra U$ in $\uuu$
\end{itemize} 
These data should satisfy a compatibility condition with respect to the $\tau$'s of section \ref{fibgradcat}.

The morphism $\phi$ is an \emph{equivalence of fibered categories} (in the 2-categorical sense) if every $\phi_U: \fff(U) \lra \GGG(U)$ is an equivalence of categories. It will be called a \emph{weak equivalence of fibered categories} \cite[Def. 2.3]{moerdijk} if every $\phi_U$ is fully faithful and \emph{locally surjective on objects}, i.e. for every $G \in \GGG(U)$ there is a covering $U_i \lra U$ and objects $F_i \in \fff(U_i)$ with $\phi_{U_i}(F_i) \cong G|_{U_i}$.

We will denote the 1-category of fibered categories and morphisms between them by $\mathsf{Fib}$ and we will denote the full subcategories of prestacks and of stacks by $\mathsf{Prestack}$ and $\mathsf{Stack}$ respectively.
For a fibered category $\fff$, an \emph{associated prestack} $ap(\fff)$ (resp. an \emph{associated stack} $as(\fff)$) is by definition a reflection of $\fff$ in $\mathsf{Prestack}$ (resp. in $\mathsf{Stack}$). Both reflections exist for a fibered category $\fff$ (see \cite{moerdijk}).
If $\ooo$ is a presheaf of rings on $\uuu$, then $ap(\fff)$ is its sheafication. 
If $\ooo$ is a sheaf of rings on $\uuu$, then $as(\fff)$ is given by \begin{equation}
\label{verstacking}
U\mapsto \{\text{locally free rank one $\ooo_U$-modules}\}
\end{equation}

\subsection{Sheaves on an $R$-linear fibered category}\label{parsheaves}

Let $\uuu \subset \mathrm{Open}(X)$ be a basis of $X$.
Our principal model of a stack of abelian categories is the stack $\rMod = \rMod(\underline{R})$ of sheaves of $R$-modules on $\uuu$. The stack $\rMod$ is defined as follows:
\begin{itemize}
\item $\rMod(U)$ is the category $\Mod(\underline{R}|_U)$ of sheaves of $R$-modules on $\uuu/U$
\item for $V \subset U$, $\rMod(U) \lra \rMod(V): F \longmapsto F|_V$ is given by restriction, i.e. $F|_V(W) = F(W)$
\item all the natural isomorphisms $\tau$ are identities
\end{itemize}
We have $\rMod|_U = \rMod(\underline{R}|_U)$.

Let $\aaa$ be an $R$-linear fibered category on $\uuu$. A \emph{sheaf} on $\aaa$ is a morphism of fibered categories $F: \aaa \lra \rMod$. Sheaves on $\aaa$ constitute an abelian category which we denote $\Mod(\aaa) = \Fib(\aaa, \rMod)$. The stack $\rMod(\aaa) = \fff ib(\aaa, \rMod)$ of sheaves on $\aaa$ is defined as follows:
\begin{itemize}
\item $\rMod(\aaa)(U) = \Mod(\aaa|_U) = \Fib(\aaa|_U, \rMod|_U)$
\item for $i:V \lra U$, we have an obvious restriction functor $i^{\ast}: \Mod(\aaa|_U) \lra \Mod(\aaa|_V)$ and we have obvious $\tau$'s
\end{itemize}
An object $A_U \in \aaa(U)$ defines a sheaf $\rHom(-,A_U): \aaa|_U \lra \rMod|_U: B_V \longmapsto \rHom(B_V,A_U|_V)$.
We have a morphism of fibered categories
$$\aaa \lra \rMod(\aaa): A_U \longmapsto \rHom(-,A_U)$$
and we will often abusively denote $\rHom(-,A_U)$ simply by $A_U$.

\begin{proposition}\label{assocmod}
Let $\aaa$ be a fibered category and let $ap(\aaa)$ and $as(\aaa)$ denote the associated prestack and stack respectively. We have equivalences of categories $\Mod(\aaa) \cong \Mod(ap(\aaa)) \cong \Mod(as(\aaa))$ and equivalences of stacks ${\rMod}(\aaa) \cong {\rMod}(ap(\aaa)) \cong {\rMod}(as(\aaa))$.
\end{proposition}

\begin{proof}
Immediate from the definition of $\Mod(\aaa)$ since $\rMod(\underline{R})$ is a stack.
\end{proof}

\subsection{The stack of sheaves in terms of localizing subcategories}\label{itoloc}

Let us first fix some notation and terminology. Let $\sss$ be a full subcategory in an arbitrary abelian category $\ccc$. The category $\sss ^{\perp}$ is by definition the full subcategory of $\ccc$ with $\Ob(\sss ^{\perp}) = \{ C \in \ccc \,| \, \Hom(\sss, C) = 0 = \Ext^1(\sss, C)\}$ whereas $^{\perp}\sss$ has $\Ob(^{\perp}\sss) = \{ C \in \ccc \, | \, \Hom(C, \sss) = 0 = \Ext^1(C, \sss)\}$. Recall that $\sss$ is called \emph{Serre} if it is closed under subquotients and extensions. In this case the quotient $q: \ccc \lra \ccc/\sss$ exists. $\sss$ is called \emph{localizing} if $q$ has a right adjoint $i$. If this is the case, $i$ is necessarily fully faithful and yields an equivalence $\ccc/\sss \cong \sss ^{\perp}$. In a Grothendieck category $\ccc$, a Serre subcategory is localizing if and only if it is closed under coproducts.

Let $\uuu \subset \mathrm{Open}(X)$ be a basis of $X$.
We will now turn to some more specific aspects of the stack $\rMod(\aaa)$ and the abelian category $\Mod(\aaa)$ over a fibered category $\aaa$ on $\uuu$. The situation is a copy of the situation for  $\rMod(\underline{R})$ and $\Mod(\underline{R})$. For $V \in \uuu$, we have an inclusion $j: \uuu/V \lra \uuu$. Consider the restriction functor $j^{\ast}: \Mod(\aaa) \lra \Mod(\aaa|_V)$, its fully faithful right adjoint $j_{\ast}: \Mod(\aaa|_V) \lra \Mod(\aaa)$ with $$(j_{\ast}F)(A_W)(W') = \lim_{U \in \uuu, U \subset V \cap W'}F(A_W|_U)(U)$$ and its fully faithful exact left adjoint $j_{!}:\Mod(\aaa|_V) \lra \Mod(\aaa)$ for which $(j_{!}F)(A_W)$ is the sheafication of $(j^p_{!}F)(A_W)$ with
$$(j^p_{!}F)(A_W)(W') = \begin{cases} F({A_W|}_{W'})(W') & \text{if}\, \, \, W' \subset V \\ 0 & \text{else} \end{cases}$$
If we take $\aaa$ to be $\aaa|_U$ on $\uuu/U$, an inclusion $j: V \subset U$ corresponds to $j: \uuu/V \lra \uuu/U$ and to functors $j^{\ast}: \Mod(\aaa|_U) \lra \Mod(\aaa|_V)$, $j_{\ast}, j_! : \Mod(\aaa|_V) \lra \Mod(\aaa|_U)$ as above.

In $\ccc = \Mod(\aaa)$ consider the subcategories
$$\zzz = \zzz(V) = \Kern(j^{\ast})$$
$$\sss = \sss(V) = j_!(\Mod(\aaa|_V))$$
$$\LLL = \LLL(V) = j_{\ast}(\Mod(\aaa|_V))$$
The situation can be summarized in the following way:
\begin{enumerate}
\item $\sss$ is a localizing Serre subcategory with $\sss^{\perp} = \zzz$
\item $\zzz$ is a localizing Serre subcategory with $\zzz^{\perp} = \LLL$
\item $\zzz^{\op}$ is a localizing Serre subcategory in $\ccc^{\op}$ with $^{\perp}\zzz = \sss = \{ C \, |\, \ccc(C,\zzz) = 0\}$ in $\ccc$
\item there is an equivalence of categories $\sss \cong \LLL$ compatible with the localization functors of (2) and (3)
\item the localization functor $\ccc \lra \LLL$ has an exact left adjoint $\LLL \cong \sss \subset \ccc$
\end{enumerate}
By the following Proposition, (3), (4) and (5) follow automatically from (1) and (2).
\begin{proposition}
Let $\ccc$ be a cocomplete abelian category with subcategories $\zzz$, $\sss$, $\LLL$ as in (1), (2) above. Then (3),(4) and (5) hold too.
\end{proposition}
\begin{proof}
Let us prove first that $^{\perp}\zzz = \sss = \{ C \, |\, \ccc(C,\zzz) = 0\}$. Obviously $\sss \subset {^\perp(\sss^{\perp})} = ^{\perp}\zzz$. Let $a: \ccc \lra \zzz$ be a localization functor with $\Kern(a) =\sss$ and consider $C$ with $\ccc(C,\sss^{\perp}) = 0$. There is an exact sequence $0 \lra S \lra C \lra aC$ with $S \in \sss$. By assumption, $C \lra aC$ is zero so $C = S \in \sss$. It now easily follows from \cite[Thm 4.5]{popescu} that $\zzz^{\op}$ is localizing in $\ccc^{\op}$. The quotient category $\ccc^{op}/\zzz^{\op}$ is thus equivalent to $(\zzz^{\op})^{\perp}$ in $\ccc^{op}$, and taking opposites we get that the corresponding functor $\ccc \lra ^{\perp}\zzz$ is equivalent to the quotient $\ccc \lra \ccc/ \zzz$ and to $\ccc \lra \zzz^{\perp} = \LLL$.
\end{proof}

Next we need to say a word on compatibility (see for example \cite{vanoystaeyenverschoren, verschoren1}). Consider two localizing Serre subcategories $\sss_U$ and $\sss_V$ in an abelian category $\ccc$. Put $\sss_U \ast \sss_V= \{ C \in \ccc \, |\, \exists\, S_U \in \sss_U, \, S_V \in \sss_V , \, 0 \lra S_U \lra C \lra S_V \lra 0 \, \}$. $\sss_U$ and $\sss_V$ are called \emph{compatible} if $\sss_U \ast \sss_V= \sss_V \ast \sss_U$. In this event, this expression is the smallest localizing Serre subcategory containing $\sss_U$ and $\sss_V$. 
Let $a_U: \ccc \lra \sss_U^{\perp}$ and $a_V: \ccc \lra \sss_V^{\perp}$ be the corresponding localization functors with right adjoint inclusion functors $i_U, i_V$ and $q_U = i_Ua_U, q_V = i_Va_V$. We will use the following incarnations of compatibility:
\begin{proposition} \cite{verschoren1, verschoren2} \label{compprop}The following are equivalent:
\begin{enumerate}
\item $\sss_U$ and $\sss_V$ are compatible
\item for all $S_U \in \sss_U$ we have $a_V(S_U) \in \sss_U$ and for all $S_V \in \sss_V$ we have $a_U(S_V) \in \sss_V$
\item $q_Uq_V = q_Vq_U$
\end{enumerate}
In this event:
\begin{enumerate} \item $q_Uq_V$ defines a localization with $\Kern(q_Uq_V) = \sss_U \ast \sss_V$ (and $(\sss_U \ast \sss_V)^{\perp} = \sss_U^{\perp} \cap \sss_V^{\perp}$)
\item $a_V$ can be restricted to a functor $a_V|_U: \sss_U^{\perp} \lra (\sss_U \ast \sss_V)^{\perp}$, which is the localization functor left adjoint to inclusion
$$\xymatrix{{\ccc} \ar[r]^{a_V} & {\sss_V^{\perp}} \\ {\sss_U^{\perp}} \ar[r]_-{a_V|_U} \ar[u] & {(\sss_U \ast \sss_V)^{\perp}} \ar[u]}$$
\item if $\sss_U^{\op}$ and $\sss_V^{\op}$ are localizing in $\ccc^{\op}$, then they are compatible
\end{enumerate}
\end{proposition}

Suppose $\uuu$ is closed under intersections. In $\Mod(\aaa)$, consider the localizing Serre subcategories $\zzz(U), \zzz(V)$ and $\zzz(U \cap V)$ for open subsets $U, V \in \uuu$. We are in the situation that $\zzz(U)$ and $\zzz(V)$ are compatible with $\zzz(U) \ast \zzz(V) =  \zzz(U \cap V)$. The square above takes the familiar form
$$\xymatrix{{\Mod(\aaa)} \ar[r]^{i_{V}^{\ast}} & {\Mod(\aaa|_V)}\\ {\Mod(\aaa|_U)} \ar[u]^{i_{U,\ast}} \ar[r]_{i_{U,U\cap V}^{\ast}} & {\Mod(\aaa|_{U\cap V})} \ar[u]_{i_{V,U\cap V,\ast}}}$$
Similarly, the square corresponding to the compatibility of $\zzz(U)^{\op}$ and $\zzz(V)^{\op}$ takes the form
$$\xymatrix{{\Mod(\aaa|)} \ar[r]^{i_{V}^{\ast}} & {\Mod(\aaa|_V)}\\ {\Mod(\aaa|_U)} \ar[u]^{i_{U,!}} \ar[r]_{i_{U,U\cap V}^{\ast}} & {\Mod(\aaa|_{U\cap V})} \ar[u]_{i_{V,U\cap V,!}}}$$
We end this section with two definitions.

\begin{definition}\label{thedefinitions}
Let $\uuu \subset \mathrm{Open}(X)$ be a full subcategory which is closed under intersections and let $\ccc$ be a fibered category on $\uuu$.
\begin{enumerate}
\item $\ccc$ is called a \emph{fibered category of localizations} if
\begin{itemize}
\item the categories $\ccc(U)$ are Grothendieck
\item for every $U \subset W$,  $\zzz_W(U) = \Kern(j^{\ast}_{W,U})$ is a localizing Serre subcategory in $\ccc(W)$
\item for every $U, V \subset W$, $\zzz_W(U)$ and $\zzz_W(V)$ are compatible and $\zzz_W(U \cap V) = \zzz_W(U) \ast \zzz_W(V)$
\end{itemize}
\item a fibered category of localizations $\ccc$ is called \emph{complemented} if for every $U \subset W$ there is a localizing Serre subcategory $\sss_W(U)$ with $\sss_W(U)^{\perp} = \zzz_W(U)$
\end{enumerate}
\end{definition}

For a fibered category of localizations $\ccc$ and $i: U \subset W$ we denote the right adjoint of the restriction functor $i^{\ast}_{W,U}: \ccc(W) \lra \ccc(U)$ by $i_{W,U,\ast}: \ccc(U) \lra \ccc(W)$. If $\ccc$ is complemented, we denote the exact left adjoint of $i^{\ast}_{W,U}$ by $i_{W,U,!}: \ccc(U) \lra \ccc(W)$.

\subsection{Some notions of flatness}\label{parflat}

Let $\uuu$ be an arbitrary base category and let $\AAA$ be an $R$-linear $\uuu$-graded category. We say that $\AAA$ is \emph{flat} if all occuring $\Hom$-modules $\AAA_f(A,A')$ are flat.
If $\aaa$ is an $R$-linear fibered category on $\uuu$, then $\aaa$ is flat as a graded category if and only if all the modules $\aaa(U)(A,A')$ are flat. 

Let $\uuu$ be a basis of $X$ and let $\aaa$ be a prestack on $\uuu$. We say that $\aaa$ is \emph{locally flat} if for all $A, A' \in \aaa(U)$ the sheaf $\rHom(A,A')$ is flat as an object of the $R$-linear abelian category $\Mod(\underline{R}|_U)$. Obviously if $\aaa$ is flat, then $\aaa$ is locally flat.
Let $\ccc$ be an $R$-linear fibered category consisting of abelian categories $\ccc(U)$. We say that $\ccc$ is \emph{flat abelian} if all the categories $\ccc(U)$ are flat abelian categories.

\begin{proposition}\label{platjes}
Let $\aaa$ be an $R$-linear prestack on a basis $\uuu$ of $X$. The following are equivalent:
\begin{enumerate}
\item $\aaa$ is a locally flat prestack
\item for every $A \in \aaa(U)$, $\rHom(-,A)$ is a flat object of $\Mod(\aaa|_U)$
\item $\rMod(\aaa)$ is a flat abelian stack
\item $\Mod(\aaa)$ is a flat abelian category
\end{enumerate}
\end{proposition}  

\begin{proof}
For $F \in \Mod(\aaa|_U)$ or $F \in \Mod(\aaa)$, $X \in \mmod(R)$, $A' \in \aaa(V)$ we have $(X \otimes_R F)(A') = X \otimes_R F(A')$ where the right hand side is computed in $\Mod(V)$. Consequently, $(X \otimes_R \rHom(-,A))(A') = X \otimes_R \rHom(A', A|_V)$ computed in $\Mod(V)$. This proves the equivalence of (1) and (2). 
To prove the equivalence of (2), (3) and (4), we first note that $\Hom_R(X,F)(A')(W) = \Hom_R(X, F(A')(W))$ because of the way limits are computed in $\Mod(V)$. Consequently, $F$ is coflat if and only if every $F(A')(W)$ is coflat. 
Now $F(A')(W) = \Mod(\aaa|_W)(\rHom(-,A'|_W), F|_W)$, so 
\begin{equation}
 \Mod(R)(X, F(A')(W)) = \Mod(\aaa|_W)(X \otimes_R \rHom(-,A'|_W), F|_W).
 \end{equation}
 Suppose (2) holds and $F$ is injective in $\Mod(\aaa)$. Then because of the existence of $j_!: \Mod(\aaa|_W) \lra \Mod(\aaa)$, $F|_W$ is injective in $\Mod(\aaa|_W)$, so coflatness of $F$ follows from flatness of $\rHom(-, A'|_W)$ and we arrive at (4). Now if $\Mod(\aaa)$ is flat then all its localizations $\Mod(\aaa_W)$ are flat as well and we arrive at (3). To see that (3) implies (2), consider $A \in \aaa(U)$. It suffices to prove that $\Mod(\aaa|_U)(X \otimes \rHom(-, A), E)$ is exact in $X$ for every injective $E \in \Mod(\aaa|_U)$. Since $E$ is coflat by assumption, the result follows once again by the equation (1) above.
\end{proof}

\subsection{The conditions (G), (F) and (FF) for morphisms of fibered categories}\label{parGFFF}
Covering systems were introduced in \cite{lowen1} as generalisations of the ``pretopologies'' of \cite{artingrothendieckverdier1}.
Let $\ccc$ be an arbitrary category. A \emph{covering system} $\ttt$ on $\ccc$ is given by specifying for every $C$ in $\ccc$ a collection $\ttt(C)$ of coverings of $C$. A \emph{covering} is by definition a collection of maps $C_i \lra C$ in $\ccc$. These coverings have to satisfy the following transitivity property: if $(C_i \lra C)_i$ is a covering of $C$ and $(C_{ij} \lra C_i)_j$ are coverings of $C_i$, then the collection of compositions $(C_{ij} \lra C_i \lra C)_{ij}$ is a covering of $C$. Also, every single identity morphism $1_C: C \lra C$ has to be a covering. This last requirement was not included in the definition of \cite{lowen1}, but is added here for convenience. Covering systems can be used both on additive and on non-additive categories, and in fact we will use both in this paper. If the underlying category $\ccc$ is additive, the notions of Grothendieck topology, site, sheaves etc. are all (implicitly) replaced by their additive versions.

A covering system $\ttt$ on a category defines a Grothendieck topology as described in \cite[Theorem 4.2]{lowen1}, i.e. for a subfunctor of $\ccc(-,C)$ to be a covering, every pullback of the subfunctor along a map $D \lra C$ has to contain a $\ttt$-covering of $D$.  We will say that a covering system is a topology if the collection of all the subfunctors generated by $\ttt$-coverings is a Grothendieck topology. We will say that two covering systems $\ttt_1$ and $\ttt_2$ are equivalent if they define the same Grothendieck topology. A category with a covering system or a Grothendieck topology $(\uuu, \ttt)$ will be called a site. The category $\Sh(\uuu, \ttt)$ of sheaves over the site is by definition the category of sheaves for the Grothendieck topology associated to the covering system. 

Covering systems can easily be induced along a functor in both directions. 
If $u: \uuu \lra \vvv$ is a functor, a covering system $\ttt$ on $\vvv$ yields an induced covering system $\ttt_u$ on $\uuu$: a collection $(U_i \lra U)$ is covering for $\ttt_u$ if and only if the collection $(u(U_i) \lra u(U))$ is covering for $\ttt$. Conversely, a covering system $\ttt$ on $\uuu$ yields an image covering system $u(\ttt)$ on $\vvv$ containing precisely the images $(u(U_i) \lra u(U))$ of $\ttt$-coverings $(U_i \lra U)$.

Of course, prestacks, stacks and sheaves over them can be defined with respect to a base site $(\uuu, \ttt)$ instead of a basis $\uuu$ of a topological space (which we implicitly endow with the natural covering system inherited from $\mathrm{Open}(X)$, for which $U_i \lra U$ is covering if and only if $\cup_iU_i = U$). 
We will not consider this more general setting in this paper when it comes to our base category $\uuu$, but we will now describe how we can ``lift'' the natural covering system of $\uuu$ to any fibered category over $\uuu$. In this way we naturally encounter sites which are no longer bases of topological spaces. 

Suppose $\ttt$ is a covering system on $\uuu$ and $\aaa$ is a fibered $\uuu$-graded category (with $F: \Ob(\aaa) \lra \Ob(\uuu)$). There is an induced covering system (denoted by $\ttt_{\aaa}$ or simply $\ttt$) on the associated additive category $\AAA$. A covering of an object $A$ is by definition a collection of cartesian morphisms $a_i \in \AAA_{u_i}(A_i, A)$ for which the collection $u_i: F(A_i) \lra F(A)$ is a $\ttt$-covering. It is readily seen from the definition of cartesian morphisms that if $\ttt$ is a topology on $\uuu$, then $\ttt_{\aaa}$ is a topology on $\AAA$.
More generaly, different covering systems $\ttt_U$ on the fibers $\aaa(U)$ can be glued together to a covering system on $\AAA$ in a similar way, but we will not need this more general construction for our purpose.

The conditions (G), (F) and (FF) were introduced in \cite{lowen1} for a morphism between additive sites $u: (\UUU, \ttt) \lra (\VVV, \zzz)$. In \cite{lowen1}, the covering system $\ttt$ is always induced by $\zzz$. We will drop this assumption here. We start by giving a definition for presheaves.
Consider a morphism $\eta: P \lra P'$ between \emph{presheaves} on $(\UUU, \ttt)$. Let $a: \Mod(\UUU) \lra \Sh(\UUU,\ttt)$ be sheafication. 

\begin{enumerate}
\item[(F)] We say that $\eta$ satisfies (F) if $a(\eta)$ is epi, i.e. if the following holds: for every $y \in P'(U)$, there is a covering $U_i \lra U$ such that $y|_{U_i}$ is in the image of $\eta_{U_i}: P(U_i) \lra P'(U_i)$.
\item[(FF)] We say that $\eta$ satisfies (FF) if $a(\eta)$ is mono, i.e. if the following holds: for every $x \in P(U)$ with $\eta_U(x) = 0$, there is a covering $U_i \lra U$ with $x|_{U_i} = 0$.
\end{enumerate}
We will now formulate the conditions for $u: (\UUU, \ttt) \lra (\VVV, \zzz)$:
\begin{enumerate}
\item[(G)] We say that $u$ satisfies (G) if for every $V \in \VVV$, there is a $\zzz$-covering $u(U_i) \lra V$.
\item[(F),(FF)] We say that $u$ satisfies (F) (resp. (FF)) if for every $U \in \UUU$, the morphism of presheaves $\UUU(-,U) \lra \VVV(u(-),u(U))$ on $(\UUU, \ttt)$ satisfies (F) (resp. (FF)).
\end{enumerate}
If $u: \UUU \lra \ccc$ is an additive functor into a Grothendieck category, we say that $u$ satisfies (G), (F) and (FF) if $u: (\UUU, \ttt_u) \lra (\ccc,\ttt)$ does, where $\ttt$ is the covering system of all epimorphic families.

The following was shown in \cite{lowen1}:
\begin{theorem}\label{gengabrpop}
Consider an additive functor $u: \UUU \lra \ccc$ from a small pre-additive category to a Grothendieck category. Let $\ttt$ be the covering system on $\ccc$ consisting of all epimorphic families of morphisms. The following are equivalent:
\begin{enumerate}
\item $u$ satisfies the conditions (G), (F) and (FF).
\item $\ttt_u$ is a topology on $\UUU$ yielding an equivalence of categories $\Sh(\UUU, \ttt_u) \cong \ccc$.
\end{enumerate}
\end{theorem}

Next we will give an interpretation of conditions (G), (F) and (FF) for a morphism of fibered categories on a fixed site. Consider a morphism of $R$-linear fibered categories $p: \bbb \lra \aaa$ on a site $(\uuu,\ttt)$, and let $\AAA$ and $\BBB$ be the associated additive categories of $\aaa$ and $\bbb$.

\begin{proposition}\label{GFFFfib}
The following are equivalent:
\begin{enumerate}
\item the morphism $p: (\BBB, \ttt_{\bbb}) \lra (\AAA, \ttt_{\aaa})$ satisfies (F) (resp. (FF))
\item for every $B, B' \in \bbb(U)$, the morphism of presheaves $\rHom_{\bbb}(B,B') \lra \rHom_{\aaa}(p(B),p(B'))$ on $(\uuu|_U,\ttt|_U)$ satisfies (F) (resp. (FF))\qed
\end{enumerate}
\end{proposition}

If the equivalent conditions of Proposition \ref{GFFFfib} hold, we will say that $p: \bbb \lra \aaa$ satisfies (F) (resp (FF)) (with respect to $\ttt$ on $\uuu$).

\begin{proposition}
The following are equivalent:
\begin{enumerate}
\item the morphism $p: (\BBB, \ttt_{\bbb}) \lra (\AAA, \ttt_{\aaa})$ satisfies (G)
\item the morphism $p$ is locally surjective \qed
\end{enumerate}
\end{proposition}

\begin{proposition}\label{GFFFstackequiv}
If $p: \bbb \lra \aaa$ satisfies (G), (F) and (FF), then the associated morphism between stacks $as(p): as(\bbb) \lra as(\aaa)$ is an equivalence.
If $\aaa$ is a prestack, the following are equivalent:
\begin{enumerate}
\item $p: \bbb \lra \aaa$ satisfies (G), (F) and (FF)
\item $ap(p): ap(\bbb) \lra \aaa$ is a weak equivalence
\item $as(p): as(\bbb) \lra as(\aaa)$ is an equivalence\qed
\end{enumerate}
\end{proposition}

We end this section with a notion which will be used in the next section \S \ref{algbfibcats}. We  say that $p$ \emph{locally reflects isomorphisms} if for every morphism $b: B \lra B' \in \bbb(U)$ with $p(b)$ an isomorphism in $\aaa(U)$, there is a covering $U_i \lra U$ with $b|_{U_i}$ an isomorphism for every $i$.

\begin{proposition}\label{FFFlociso}
If $p: \bbb \lra \aaa$ satisfies (F) and (FF), then $p$ locally reflects isomorphisms.
\end{proposition}

\begin{proof}
Consider $b: B \lra B' \in \bbb(U)$ with $p(b)$ an isomorphism in $\aaa(U)$. Let $a: p(B') \lra p(B)$ be an inverse isomorphism to $p(b)$. By (F) there is a covering $V \lra U$ on which we have morphisms $c_V: B'|_V \lra B|_V$
for which $p(b|_V)$ and $p(c_V)$ are inverse isomorphisms. Consequently, by (FF), there are further coverings $W \lra V$ on which $b|_V|_W$ and $c_V|_W$ are inverse isomorphisms.
\end{proof}

\subsection{Algebroid fibered categories}\label{algbfibcats}

Algebroids and algebroid fibered categories were introduced by Kontsevich in \cite{kontsevich1}.

\begin{definition}
An $R$-linear category $\AAA$ is called an \emph{algebroid} if it is nonempty and all its objects are isomorphic.
An $R$-linear fibered category $\aaa$ on a site $(\uuu, \ttt)$ is called \emph{algebroid} (or an \emph{$R$-algebroid fibered category}) if the following hold:
\begin{enumerate}
\item For $U \in \uuu$, there is a covering $U_i \lra U$ with $\aaa(U_i) \neq \varnothing$.
\item For $A, A' \in \aaa(U)$, there is a covering $U_i \lra U$ with $A|_{U_i} \cong A'|_{U_i}$.
\end{enumerate}
\end{definition}

The fibered category associated to a prescheaf $\ooo$ is obviously algebroid, all the one-object categories $\ooo(U)$ being algebroids. One can consider an intermediate notion between presheaves and algebroid fibered categories, which involves the 2-category $2\Alg(R)$ of $R$-algebras. This 2-category is the ``2-full'' subcategory of $2\Alg(R)$ with as objects ($0$-cells) $R$-algebras (considered as one-object categories). Explicitly, a $1$-cell between $R$-algebras is just an $R$-algebra morphism, and a $2$-cell $\eta: f \lra g$ between $1$-cells $f, g: A \lra B$ is an element $\eta \in B$ such that for all $a \in A$ we have $g(a)\eta = \eta f(a)$.

\begin{definition}\label{twisted}
A \emph{twisted presheaf} of $R$-algebras on a base category $\uuu$ is a pseudofunctor $$\ooo: \uuu \lra 2\Alg(R)$$ 
\end{definition}
Explicitly, a twisted presheaf $\ooo$ consists of $R$-algebras $\ooo(U)$ for $U \in \uuu$, restriction morphisms $i^{\ast}: \ooo(U) \lra \ooo(V)$ for $i: V \lra U$, elements $\tau_{i,j} \in \ooo(W)$ for $j: W \lra V$, satisfying the following ``cocycle condition'' for $k: Z \lra W$:
$$k^{\ast}(\tau_{i,j}) \tau_{ij,k} = \tau_{j,k} \tau_{i,jk}$$

\begin{proposition}
Suppose $\aaa$ is a fibered $\uuu$-$R$-category such that all the categories $\aaa(U)$ for $U \in \uuu$ are algebroids. There exists a twisted presheaf of $R$-algebras $\ooo$ on $\uuu$ and an equivalence of fibered categories $\ooo \cong \aaa$.
\end{proposition}

\begin{proof} Consider $\aaa$ as a graded category. For every $U \in \uuu$, pick one object $A_U \in \aaa(U)$ and let $\ooo$ bet the full graded subcategory of $\aaa$ spanned by the objects $A_U$. By composing cartesian morphisms with isomorphisms, $\ooo$ is readily seen to be fibered too, and $\ooo \cong \aaa$.
\end{proof}

The next few propositions give some relations between algebroid fibered categories and the conditions (G), (F) and (FF) of the previous section \S \ref{parGFFF} for morphism between them. They will be used in the main section \S \ref{abdefstack} of the second chapter in the context of deformations $\bbb \lra \aaa$ of fibered categories (Proposition \ref{algliftdefo}, Proposition \ref{algebroidclue}).

\begin{proposition}\label{alglift}
Consider a morphism of $R$-linear fibered categories $p: \bbb \lra \aaa$ on a site $(\uuu,\ttt)$ which satisfies (G) and (F) and locally reflects isomorphisms. If $\aaa$ is algebroid, then so is $\bbb$.
\end{proposition}

\begin{proof}
For $U \in \uuu$, let $U_i \lra U$ be a cover with $\aaa(U_i) \neq \varnothing$. Choose objects $A_i \in \aaa(U_i)$ and let $U_{ij} \lra U_i$ be covers on which $A_i|_{U_{ij}} \cong p(B_{ij})$. In particular, the cover $U_{ij} \lra U$  is such that $\bbb(U_{ij}) \neq \varnothing$. 
Now consider $B, B' \in \bbb(U)$. For $p(B), p(B') \in \aaa(U)$, there exists a cover $U_i \lra U$ with $p(B)|_{U_i} \cong p(B')|_{U_i}$. Consequently, we have isomorphisms $a_i: p(B|_{U_i}) \lra p(B'|_{U_i})$. By (F), there are covers $U_{ij} \lra U_i$ with $a_i|_{U_{ij}} \cong p(b_{ij})$ for some $b_{ij}: B|_{U_{ij}} \lra B'|_{U_{ij}}$. By assumption, the maps $b_{ij}$ are isomorphisms on some further covering.
\end{proof}

\begin{proposition}\label{algliftGFFF}
Consider a morphism of $R$-linear fibered categories $p: \bbb \lra \aaa$ on a site $(\uuu,\ttt)$ which satisfies (G), (F) and (FF).  If $\aaa$ is algebroid, then so is $\bbb$.
\end{proposition}

\begin{proof}
Immediate from Propositions \ref{alglift} and \ref{FFFlociso}.
\end{proof}

\begin{proposition}\label{algG}
Consider a morphism $p: \bbb \lra \aaa$ of fibered categories. If $\bbb$ is algebroid, then $\aaa$ is algebroid if and only if $p$ satisfies (G).
\end{proposition}

\begin{proof}
Suppose $\aaa$ is algebroid.
Consider $A \in \aaa(U)$. Take a covering $U_i \lra U$ such that $\bbb(U_i) \neq \varnothing$. Choose $B_i \in \bbb(U_i)$. Take coverings $U_{ij} \lra U_i$ on which $p(B_i)|_{U_{ij}} \cong A_{U_{ij}}$.
Suppose $p$ satisfies (G). Take a covering $U_i \lra U$ such that $\bbb(U_i) \neq \varnothing$. Then certainly $\aaa(U_i) \neq \varnothing$. Take $A, A' \in \aaa(U)$. Take a covering on which both $A|_V \cong p(B)$ and $A'|_V \cong p(B')$ and take a further cover on which $B$ and $B'$ become isomorphic.
\end{proof}

\subsection{Characterization of the stack of sheaves on $\aaa$}\label{parkarmod}
Let $\aaa$ be a fibered category on a topological space $X$. In this section we characterize the stack $\rMod(\aaa)$ in terms of the inclusion $\aaa \lra \rMod(\aaa)$, using additive sheaf categories. In the second chapter, this characterization (which is formulated in Theorem \ref{karmodthm}) will be lifted under deformation.

We start this section with a proof of the standard characterization of module categories among Grothendiek categories. This proof will be our inspiration for Theorem \ref{karmodthm}.

\begin{theorem}\label{baby}
Consider an additive functor $\AAA \lra \ccc$ from a small pre-additive category $\AAA$ to a Grothendieck category $\ccc$. Let $\ttt_{\mathrm{epi}}$ be the covering system on $\AAA$ induced by the covering system of all epimorphic families on $\ccc$, and let $\ttt_{\mathrm{triv}}$ be the trivial covering system on $\AAA$, for which the only coverings are identities. The following are equivalent:
\begin{enumerate}
\item $\AAA \lra \ccc$ yields an equivalence $\Mod(\AAA) \cong \ccc$
\item 
\begin{enumerate}\item$\AAA \lra \ccc$ satisfies (G), (F) and (FF)
\item there is an equivalence of covering systems $\ttt_{\mathrm{epi}} \cong \ttt_{\mathrm{triv}}$
\end{enumerate}
\item $\AAA \lra \ccc$ is fully faithful and the (images of) objects of $\AAA$ form a set of finitely presented projective generators in $\ccc$
\end{enumerate}
\end{theorem}

\begin{proof}
The equivalence of (1) and (2) immediately follows from Proposition \ref{gengabrpop} and the implication (1) implies (3) is obvious. To prove that (3) implies (2), we consider an arbitrary covering $f_i: A_i \lra A$ for $\ttt_{\mathrm{epi}}$. By definition, we have a $\ccc$-epimorphism $\coprod_iA_i \lra A$. By (2), this epimorphism splits through a finite sub-coproduct $A \lra \coprod_{i \in J}A_i$. Consequently, since $\AAA \lra \ccc$ is fully faithful, we get morphisms $g_i: A \lra A_i$ for $i \in J$ with $\sum_{i \in J}f_ig_i = 1_A$. This finishes the proof.
\end{proof}

As a first step, we describe $\Mod(\aaa)$ as an additive sheaf category:
\begin{proposition}\label{modsh}
Let $\aaa$ be a fibered category on a basis $\uuu$ of $X$, let $\AAA$ be the additive category of $\aaa$ and let $\ttt_{\aaa}$ be the topology on $\AAA$ induced by the standard topology $\ttt$ on $\uuu$ (see \S \ref{parGFFF}). There is an equivalence of categories $$\Mod(\aaa) \cong \Sh(\AAA, \ttt_{\aaa})$$
\end{proposition}

\begin{proof}
We can easily give functors $\varphi: \Mod(\aaa) \lra \Sh(\AAA, \ttt_{\aaa})$ and $\psi: \Sh(\AAA, \ttt_{\aaa}) \lra \Mod(\aaa)$ constituting an equivalence. Let $F: \aaa \lra \rMod(\uuu)$ be a sheaf on $\aaa$. Then we define $\varphi F(A_U) = F(A_U)(U)$. For a morphism $A_W \lra A_V|_{V_i}$ from $A_W$ to $A_V$ in $\AAA$, we get a morphism $F(A)(V) \lra F(A)(W) \cong F(A|_{W})(W) \lra F(A_W)(W)$ as required. Conversely, for an additive sheaf $G: \AAA \lra \Ab$, we put $\psi G(A_V)(W) = G(A_V|_W)$. A map $A_V \lra A_V'$ in $\aaa(V)$ yields morphisms $G(A_V|_W) \lra G(A_V'|_W)$ as required. It is easily seen that by the definition of $\ttt_{\aaa}$, sheaves are mapped to sheaves by both $\varphi$ and $\psi$, and that they are inverse equivalences.
\end{proof}

Let $\ccc$ be a \emph{complemented stack of localizations} (see Definition \ref{thedefinitions}) on $\uuu = \mathrm{Open}(X)$ and consider a morphism $\aaa \lra \ccc$. Let $\AAA_U$ be the additive category associated to the fibered graded category $\aaa|_U$ (see \S \ref{fibgradcat}).  For every $U \in \uuu$, there is an induced $$\AAA_U \lra \ccc(U): A_V \longmapsto j_{U,V,!}{A_V}$$ with the following prescription for morphisms if $W \subset V$:
$$\AAA_U(A_W,A_V) = \aaa(W)(A_W, A_V|_W) \lra  \ccc(W)(A_W, A_V|_W) \cong \ccc(U)(j_{!}A_W, j_!A_V)$$

The category $\AAA_U$ naturally carries two topologies:
\begin{enumerate}
\item the topology $\ttt_U = \ttt_{\aaa|_U}$ which is induced by the standard covering system $\ttt$ on $\mathrm{Open}(X)$
\item the topology $\ttt_{\mathrm{epi},U}$ induced by the inclusion $\AAA_U \lra \ccc(U)$, where $\ccc(U)$ is endowed with the covering system $\ttt_{\mathrm{epi}}$ of all epimorphic families
\end{enumerate}

\begin{proposition}\label{inc}\label{cover}
We have an inclusion of topologies on $\AAA_U$:
$$\ttt_U \subset \ttt_{\mathrm{epi},U}$$
\end{proposition}

\begin{proof}
Consider a $\ttt_V$-covering of an object $A \in \aaa(V)$ for $V \subset U$, i.e. a collection of $\ccc(V_i)$-isomorphisms $A_i \lra A|_{V_i}$ for a covering $V_i \lra V$. To show that these morphisms are covering for $\ttt_{\mathrm{epi},U}$, it suffices that their images $j_{V,V_i!}A_i \lra A$ are epimorphic in $\ccc(V)$. Suppose all compositions $j_{V,V_i!}A_i \lra A \lra M$ are zero for $A \lra M$ in $\ccc(V)$. Equivalently, all the restrictions $A|_{V_i} \lra M|_{V_i}$ are zero, hence it follows that $A \lra M$ is zero since $\ccc$ is a stack.
\end{proof}
 
The equivalence of (1) and (2) in Theorem \ref{karmodthm} immediately follows from Proposition \ref{modsh} and Theorem \ref{gengabrpop}, and in fact this part of the theorem is sufficient for our purpose, namely for the proof of Theorem \ref{theorem} in the second chapter. However, it is possible to make the analogy with Theorem \ref{baby} complete by defining local versions of projective and finitely presented objects in a stack. These notions will be used again later on in \S \ref{abdefstack}. The proof of Theorem \ref{karmodthm} is essentially a local version of the proof of Theorem \ref{baby}, but we include it for completeness.

\begin{definition}\label{locep}
Let $A \in \ccc(U)$ be an object of a stack of ($R$-linear) abelian categories $\ccc$.
\begin{enumerate}
\item We say that $A$ is \emph{locally projective} if $\rExt^1(A,-): \ccc(U) \lra \Mod(U)$ is zero. 
\item We say that $A$ is \emph{locally finitely presented} if $\rHom(A,-): \ccc(U) \lra \Mod(U)$ preserves filtered colimits.\end{enumerate}
\end{definition}

\begin{remark}
Let $A \in \ccc(U)$ be as above.
\begin{enumerate}
\item $A$ is locally projective if and only if for every $C \in \ccc(U)$, $0 \lra \Ext_{\ccc}(A,C)$ satisfies (F). In other words, if every epimorphism $f: B|_V \lra A|_V$ in $\ccc(V)$ (with $B \in \ccc(U)$) is locally split, i.e. there is a covering $V_i \lra V$ such that for every $i$ there is a splitting $g: A|_{V_i} \lra B|_{V_i}$ with $f|_{V_i}g= 1$.
\item $A$ is locally finitely presented if and only if, for every filtered colimit  $\colim_iB_i$ in $\ccc(U)$, $\colim_i\Hom_{\ccc}(A,B_i) \lra \Hom_{\ccc}(A,\colim_iB_i)$ satisfies (F) and (FF). In other words,
\begin{enumerate}
\item  if $f: A|_V \lra \colim_i(B_i)|_V$ a morphism in $\ccc(V)$, then there is a covering $V_j \lra V$ such that for every $j$, $f|_{V_j}$ factors through $A|_{V_j} \lra (B_i)|_{V_j}$ for some $i$,
\item if $f: A|_V \lra (B_i)|_V$ is such that the composition $A|_V \lra \colim_i(B_i)|_V$ is zero, then there is a covering $V_j \lra V$ with $f|_{V_j} = 0$ for every $j$.
\end{enumerate} 
\end{enumerate}
\end{remark}

\begin{theorem}\label{karmodthm}
Let $\ccc$ be a complemented stack of localizations on $X$, $\aaa$ a fibered category on $X$, and $\aaa \lra \ccc$ a morphism of fibered categories.
The following are equivalent:
\begin{enumerate}
\item $\aaa \lra \ccc$ yields an equivalence of stacks $\rMod(\aaa) \cong \ccc$.
\item 
\begin{enumerate}
\item the maps $\AAA_U  \lra \ccc(U)$ statisfy (G), (F) and (FF)
\item there are equivalence of topologies $\ttt_{\mathrm{epi},U} \cong \ttt_U$
\end{enumerate}
\item \begin{enumerate}
\item the maps $(\AAA_U,\ttt_U)  \lra (\ccc(U), \ttt_{\mathrm{epi}})$ statisfy (G), (F) and (FF)
\item the (images of) objects of $\aaa$ are locally finitely presented and locally projective in $\ccc$
\end{enumerate}
\end{enumerate}
\end{theorem}

\begin{remark}
Note that in Theorem \ref{baby}(3), $\AAA \lra \ccc$ is fully faithful is if and only if $(\AAA, \ttt_{\mathrm{triv}}) \lra (\ccc, \ttt_{\mathrm{epi}})$ satisfies (F) and (FF).
\end{remark}

\begin{remark}
Note that even if $\aaa \lra \ccc$ is a fully faithful embedding (in the sense that all the functors $\aaa(U) \lra \ccc(U)$ are), condition (FF) in (2)(a) is automatically fulfilled, but condition (F) is not.
\end{remark}

\begin{proof}
The equivalence of (1) and (2) immediately follows from Proposition \ref{modsh} and Theorem \ref{gengabrpop}.
In the situation of (1), (3)(a) holds because (2) holds. For (2)(b) we are to show that $\Ext^1_{{\rMod}(\aaa)}(\rHom_{\aaa}(-,A),M)$ is a zero presheaf for $A \in \aaa(U)$, $M \in \Mod(\aaa|_U)$. Consider an exact sequence $0 \lra M|_V \lra E \lra \rHom_{\aaa}(-,A|_V) \lra 0$ in $\Mod(\aaa|_V)$. Since $E(A|_V) \lra \rHom_{\aaa}(A|_V,A|_V)$ is an epimorphism of sheaves, there is a covering $V_i \lra V$ on which $1: A|_{V_i} \lra A|_{V_i}$ is in the image of $E(A|_{V_i})$. This allows us to define a splitting of $E|_{V_i} \lra \rHom_{\aaa}(-,A|_{V_i})$.
Next we are to show that if $\colim_i M_i$ is a filtered colimit in $\Mod(\aaa|_U)$, then the canonical map from $\colim_i \Hom_{\Mod(\aaa|_V)}(\rHom(-,A|_V), M_i|_V)$ to $\Hom_{\Mod(\aaa|_V)}(\rHom(-,A|_V), \colim_i M_i|_V)$ satisfies (F), (FF). This map can be rewritten as $\colim_i (M_i(A)(V)) \lra (\colim_iM_i(A))(V)$. Since this map corresponds to sheafication of the \emph{presheaf} $\colim_i (M_i(A))$, this finishes the proof of (2)(b).
Finally we show that (3) implies (2).
By Proposition \ref{inc}, it suffices to show that for a covering $A_i \lra A|_{V_i}$ of $A \in \aaa(V)$ for $\ttt_{\mathrm{epi},U}$, the generated subfunctor contains a covering for ${\ttt}_U$. So we suppose that the maps $j_{V,V_i}A_i \lra A$ are epimorphic in $\ccc(V)$, i.e. we have an epimorphism $\coprod_ij_{V,V_i!}A_i \lra A$. For an open subset $W \subset V_i$, composing $A_i \lra A|_{V_i}$ with the canonical $A_i|_{W} \lra A_i|_{W}$ yields the $\AAA_U$-morphism from $A_i|_W$ to $A$ determined by $A_i|_{W} \lra A|_{V_i}|_{W} \cong A|_W$. Hence it suffices to find for every $V_i$ a covering $W \lra V_i$ for which the maps $A_i|_{W} \lra A|_{V_i}|_{W}$ are isomorphisms.
By (2), there is a covering $W_k \lra V$ on which the epimorphism $\coprod_ij_{V,V_i!}A_i \lra A$ splits, and on which the splitting is through a finite sub coproduct $\coprod_{i \in J_k}(j_{V,V_i!}A_i)|_{W_k}$. In other words, for every $k$ and $i \in J_k$, for the map $p_{ki}:j_{W_k,V_i \cap W_k!}A_i|_{V_i \cap W_k} \lra A|_{W_k}$ there is an $s_{ki}: A|_{W_k} \lra j_{W_k,V_i \cap W_k!}A_i|_{V_i \cap W_k}$ such that $\sum_{i \in J_k} p_{ki} s_{ki} = 1$.
The map $p_{ki}$ is the image in $\ccc(W_k)$ of $\bar{p}_{ki}: A_i|_{V_i \cap W_k} \lra A|_{V_i \cap W_k}$, which is a map in the generated subfunctor.
The maps $s_{ki}$ are not in the image of $\AAA_{W_k} \lra \ccc({W_k})$, but by (1)(F), there exist coverings $W_{kl} \lra W_k$ for which the canonical $s_{kl}: j_{W_{k},W_{kl}!}A|_{W_{kl}} \lra A|_{W_k}$ are such that $s_{ki} s_{kl}$ is in the image for every $i \in J_k$. Suppose $s_{ki} s_{kl}$ is the image of $\bar{s}_{kil}$. Then for every $l$, $\sum_{i \in J_k}\bar{p}_{ki}\bar{s}_{kil}$ is mapped onto $s_{kl}$. But $s_{kl}$ is by definition the image of the canonical $\bar{s}_{kl}: A|_{W_{kl}} \lra  A|_{W_{kl}}$. Hence by (1)(FF), there is a further covering $W_{klm} \lra W_{kl}$ such that the canonical morphism $\bar{s}_{klm}:A|_{W_{klm}} \lra A|_{W_{klm}}$ can be written as
$\sum_{i \in J_k}\bar{p}_{ki}\bar{s}_{kil}\bar{s}_{klm}$, hence $\bar{s}_{klm}$ is in the subfunctor generated by the original $A_i \lra A|_{V_i}$, as required.
\end{proof}

\section{Deformations of ringed spaces}

Throughout, let $S \lra R$ be a surjective ring map (between commutative, coherent rings) with nilpotent kernel $I$ (we may and will assume that $I^2 = 0$).

The main aim of this chapter is to show for a ringed space $(X, \ooo)$ how we can describe the abelian deformations of the category $\Mod(\ooo)$ of sheaves of modules over $\ooo$ in terms of a certain type of deformations of $\ooo$. To do so, we will have to consider $\ooo$ no longer as a sheaf of algebras but as an algebroid prestack of linear categories on $X$. We will consider this prestack $\ooo$ as sitting inside the stack $\rMod(\ooo)$ of sheaves of modules over $\ooo$, which has the category $\Mod(\ooo_U)$ of sheaves of modules over $\ooo_U$ as section on $U$.
In \S \ref{abdefstack}, Theorem \ref{theorem}, for appropriate notions of deformations, we obtain equivalences between
\begin{itemize}
\item deformations of the prestack $\ooo$
\item deformations of the stack $\rMod(\ooo)$
\item deformations of the abelian category $\Mod(\ooo)$
\end{itemize}
In sections \S \ref{lindef} and \S \ref{abdef} we define what we mean by these different types of deformations, and we prove some preliminary results. Theorem \ref{theorem} is based upon a technical result which allows us to ``lift'' the conditions (G), (F) and (FF) of \S \ref{parGFFF} under deformation. This is explained in \S \ref{liftGFFF}. In the last two sections \S \ref{paracyclic} and \S \ref{parafcover} we analyze what happens if $(X,\ooo)$ is a quasi-compact, separated scheme. This allows us to prove that there is an equivalence between
\begin{itemize}
\item deformations of the abelian category $\Mod(X)$
\item deformations of the abelian category $\Qch(X)$
\end{itemize}  
The parallel result for Hochschild cohomology has been demonstrated in \cite{lowenvandenbergh2}.

\subsection{Linear deformations of fibered categories}\label{lindef}

In \S \ref{fibgradcat} we explained that there are two ways to think of a fibered category. We will now show that there are obvious notions of deformations in both cases, and that these notions are equivalent.

Let $\uuu$ be an arbitrary base category and let $\AAA$ be an $\uuu$-graded $R$-linear category. An \emph{$S$-deformation of $\AAA$} is an $\uuu$-graded $S$-linear category $\BBB$ with an equivalence $R \otimes_S \BBB \lra \AAA$. Here $R \otimes_S \BBB$ is obtained by simply tensoring all modules with $R$. Note that for an ordinary $R$-linear category (which is naturally graded over the category $\uuu_{\AAA}$ with $\Ob(\uuu_{\AAA}) = \Ob(\AAA)$ and $\uuu_{\AAA}(A,A') = \{ \ast \}$ for all $A, A' \in \AAA$) this yields the definition of a deformation of linear categories as in \cite{lowenvandenbergh1}. We immediately get the following

\begin{proposition}\label{gradlin}
Suppose $\AAA$ is a $\uuu$-graded $R$-linear category over a \emph{poset} $\uuu$, and denote by the same symbol $\AAA$ the associated $R$-linear category.
There is an equivalence between
\begin{enumerate}
\item flat deformations of $\AAA$ as a $\uuu$-graded category
\item flat deformations of $\AAA$ as a linear category
\end{enumerate}
\end{proposition}

\begin{proof}
The associated $R$-linear category $\AAA$ has $\AAA(A,A') = \AAA_{\leq}(A,A')$ if $F(A) \leq F(A')$ and $\AAA(A,A') = 0$ otherwise. The result follows since the zero modules necessarily deform to zero. \end{proof}

\begin{proposition}\label{fibbew}
Let $\BBB \lra \AAA$ be a flat deformation of graded categories. If $\AAA$ is fibered, then so is $\BBB$.
\end{proposition}

\begin{proof}
Consider $f:V \lra U$ in $\uuu$ and $B \in \BBB(U)$. Take a cartesian morphism $a: A \lra R \otimes_S B$ in $\AAA_f(A, R \otimes_S B)$ and let $\bar{a}: \bar{A} \lra B$ be a lift to $\BBB$. Consider $C \in \BBB(W)$ and $g: W \lra V$. We have a commutative diagram
$$\xymatrix{ {\BBB_g(C,\bar{A})} \ar[d] \ar[r]^{\bar{a}-}& {\BBB_{fg}(C,B)} \ar[d] \\
{\AAA_g(R \otimes_S C, A)} \ar[r]_-{{a}-} & {\AAA_{fg}(R \otimes_S C, R \otimes_S B)}}$$
in which $a- = R \otimes_S (\bar{a}-)$. By assumption, $a-$ is an isomorphism, and consequently, by flatness and nilpotency, $\bar{a}-$ is too.
\end{proof}

Let $\aaa$ be an $R$-linear fibered category on $\uuu$, viewed as a pseudofunctor. A \emph{linear $S$-deformation of $\aaa$} is an $S$-linear fibered category $\bbb$ on $\uuu$ with an equivalence $R \otimes_S \bbb \lra \aaa$. Here $R \otimes_S \bbb$ is obtained by tensoring every $\bbb(U)$ with $R$ and inducing the restriction functors. The following easily follows from Proposition \ref{fibbew}.
\begin{proposition}\label{gradfib}
Suppose $\aaa$ is a fibered category over $\uuu$. There is an equivalence between
\begin{enumerate}
\item flat deformations of $\aaa$ as a graded category
\item flat deformations of $\aaa$ as a fibered category\qed
\end{enumerate}
\end{proposition}

Let $\ooo$ be a twisted presheaf of $R$-algebras on $\uuu$. An \emph{$S$-deformation of $\ooo$} is a twisted presheaf $\ppp$ of $S$-algebras on $\uuu$ with an equivalence $R \otimes_S \ppp \lra \ooo$. Here $R \otimes_S \ppp$ is obtained by tensoring every $\ppp(U)$ with $R$ and inducing the restriction functors.

\begin{proposition}\label{twistfibgrad}
Let $\bbb \lra \aaa$ be a flat deformation of fibered categories. If $\aaa(U)$ is an algebroid, then so is $\bbb(U)$. Consequently, for a twisted presheaf of $R$-algebras $\ooo$ on $\uuu$, there is an equivalence between
\begin{enumerate}
\item flat deformations of $\ooo$ as a twisted presheaf
\item flat deformations of $\ooo$ as a fibered category
\item flat deformations of $\ooo$ as a graded category
\end{enumerate}
\end{proposition}

\begin{proof} By nilpotency, an isomorphism in $\aaa(U)$ can be lifted to an isomorphism in $\bbb(U)$.
\end{proof}

Let $\uuu$ be a basis of $X$. We will now introduce a different notion of deformation for a fibered category $\aaa$, which will be the correct one to obtain an equivalence with abelian deformations of $\Mod(\aaa)$. By proposition \ref{assocmod}, changing $\aaa$ to $ap(\aaa)$ or $as(\aaa)$ doesn't change the stack $\rMod(\aaa)$. This suggests that in order to get an equivalence with deformations of $\Mod(\aaa)$, we have to deform $\aaa$ ``up to stackification''. Since the ``pointwise'' functor
$$R \otimes_S -: \mathsf{Fib} \lra \mathsf{Fib}: \bbb \longmapsto R \otimes_S \bbb$$
with $(R \otimes_S \bbb)(U) = R \otimes_S \bbb(U)$ preserves neither prestacks nor stacks, we look at
$$\mathsf{Stack} \lra \mathsf{Stack}: \bbb \longmapsto as(R \otimes_S \bbb)$$
instead. For an $R$-linear stack $\aaa$, a \emph{linear stack $S$-deformation} of $\aaa$ is by definition an $S$-linear stack $\bbb$ together with an equivalence $as(R \otimes_S \bbb) \cong \aaa$. For an arbitrary fibered category $\aaa$, we define a \emph{weak linear deformation} of $\aaa$ to be a stack deformation of $as(\aaa)$. Since we are most interested in prestacks $\aaa$ (a sheaf of algebras $\ooo$ naturally defining a prestack but not a stack), we will elaborate this a little further for prestacks.
Let $\Sigma$ denote the class of weak equivalences (\S \ref{parstacks}) in the category $\mathsf{Prestack}$. Then $\Sigma$ contains precisely the morphisms inverted by the reflection $as: \mathsf{Prestack} \lra \mathsf{Stack}$. Consequently, $\mathsf{Stack} \cong \mathsf{Prestack}[\Sigma^{-1}]$. For prestacks $\aaa$ and $\bbb$, every morphism from $\aaa$ to $\bbb$ in $\mathsf{Prestack}[\Sigma^{-1}]$ can be represented by a ``fraction'' consisting of a morphism $\aaa \lra \ccc$ to a prestack $\ccc$ and a weak equivalence $\bbb \lra \ccc$. If both parts of the fraction are weak equivalences, we will call the resulting morphism $\aaa \lra \ccc$ a \emph{stack equivalence}.
Let $\aaa$ be an $R$-linear prestack on $\uuu$. An weak linear $S$-deformation of $\aaa$ can be represented by an $S$-linear prestack $\bbb$ with an equivalence $as(R \otimes_S as(\bbb)) \cong as(\aaa)$. Since $R \otimes_S -$ preserves weak equivalences, this corresponds to an $S$-linear prestack $\bbb$ with a stack equivalence $ap(R \otimes_S \bbb) \cong \aaa$. Two weak deformations $\bbb$ and $\bbb'$ are equivalent if there is a stack equivalence $\bbb \lra \bbb'$ inducing the stack equivalence $ap(R \otimes_S \bbb) \lra ap(R \otimes_S \bbb')$. We will be most concerned with weak deformations $\bbb$ of $\aaa$ where the stack equivalence $ap(R \otimes_S \bbb) \lra \aaa$ is in fact a weak equivalence. Equivalent characterizations of this situation are given in Proposition \ref{GFFFstackequiv}.

\begin{proposition}\label{algliftdefo}
Suppose we are in one of the following cases:
\begin{enumerate}
\item $\bbb \lra \aaa$ is a linear deformation of fibered categories
\item $\bbb \lra \aaa$ is a weak linear deformation of prestacks
\end{enumerate}
If $\aaa$ is algebroid, then so is $\bbb$.
\end{proposition}

\begin{proof}
The first case follows from Proposition \ref{alglift}. For the second case we have to consider morphisms $R \otimes_S \bbb \lra ap(R \otimes_S \bbb) \lra \ccc \longleftarrow \aaa$. Then $\ccc$ is algebroid by Proposition \ref{algG}, and $R \otimes_S \bbb$ by Proposition \ref{algliftGFFF}.
\end{proof}

\subsection{Abelian deformations of fibered categories of localizations}\label{abdef}

Let $\ccc$ be an $R$-linear fibered category of abelian categories. An \emph{abelian $S$-deformation of $\ccc$} is an $S$-linear fibered category of abelian categories $\ddd$ with an equivalence of fibered categories $\ccc \lra \ddd_R$. Here $\ddd_R$ is the $R$-linear fibered category with $\ddd_R(U) = \ddd(U)_R$ and the induced restriction functors, and $\ddd(U)_R$ is the category of $R$-objects in $\ddd(U)$, i.e. those objects annihilated by the kernel $I$ of $S \lra R$.
For $\uuu = \{ \ast, 1_{\ast} \}$, this reduces to the definition of a deformation of abelian categories of \cite{lowenvandenbergh1}.

Let $\uuu \subset \mathrm{Open}(X)$ be a base category which is closed under intersections and let $\ccc$ be an $R$-linear fibered category of localizations on $\uuu$. A \emph{localization $S$-deformation of $\ccc$} is an abelian $S$-deformation which is again a fibered category of localizations. If $\uuu$ is a basis of $X$ and $\ccc$ is a stack on $\uuu$, we are interested in deformations which are again stacks. We will call them \emph{stack deformations}, or \emph{localization stack deformations} if they are at the same time localization deformations.
The remainder of this section contains preliminary results about the liftability under deformation of several features of fibered categories. They will be used in \S \ref{abdefstack} and further on.
The following Proposition shows that if $X \in \uuu$, the localization deformations of $\ccc$ can all be induced from deformations of $\ccc(X)$.

\begin{proposition}\label{largest}
Let $\ccc$ be an $R$-linear fibered category of localizations on $\uuu$ and suppose $X \in \uuu$. There is an equivalence between 
\begin{enumerate}
\item flat abelian deformations of $\ccc(X)$
\item flat localization deformations of $\ccc$
\end{enumerate}
\end{proposition}

\begin{proof}
Consider a flat abelian deformation $\ccc(X) \lra \ddd(X)$. By \cite{lowenvandenbergh1}, there is a (up to isomorphism) unique localization $j^{\ast}_{\ddd}: \ddd(X) \lra \ddd(U)$ ``lifting'' the localization $j^{\ast}_{\ccc}: \ccc(X) \lra \ccc(U)$, and this localization is obtained by lifting $\Kern( j^{\ast}_{\ccc})$ to its generated Serre subcategory in $\ddd(X)$. So it remains to show that the remaining conditions in Definition \ref{thedefinitions} are fulfilled. This follows from Proposition \ref{complift}.
\end{proof}

\begin{proposition}\label{complemented}
Let $\ccc \lra \ddd$ be a flat localization deformation. If $\ccc$ is complemented, then so is $\ddd$.
\end{proposition}

\begin{proof}
For $j:V \subset U$, consider $j_{\ddd}^{\ast}: \ddd(U) \lra \ddd(V)$. The kernel $\zzz'$ of $j_{\ddd}^{\ast}$ is the Serre subcategory generated by $\zzz = \Kern(\ccc(U) \lra \ccc(V))$ in $\ddd$. Now let $\sss$ be the localizing Serre subcategory in $\ccc$ with $\sss^{\perp} = \zzz$ and let $\sss' = \langle \sss \rangle_{\ddd}$ be the generated Serre subcategory in $\ddd$. By proposition \ref{complift}(1), $\sss'^{\perp} = \zzz'$, which finishes the proof.
\end{proof}

\begin{proposition}\label{complift}
Consider a flat deformation of abelian categories $j: \ccc \lra \ddd$ and localizing Serre subcategories $\sss'$ and $\zzz'$ in $\ddd$ with $\sss = \sss' \cap \ccc$, $\zzz = \zzz' \cap \ccc$.
\begin{enumerate}
\item if $\sss^{\perp} = \zzz$, then $\sss'^{\perp} = \zzz'$.
\item is $\sss$ and $\zzz$ are compatible, then so are $\sss'$ and $\zzz'$ and $(\sss' \ast \zzz')\cap \ccc = \sss \ast \zzz$.
\end{enumerate}
\end{proposition}

\begin{proof}
(1) Suppose $\sss^{\perp} = \zzz$.  We start with the inclusion $\zzz' \subset \sss'^{\perp}$. Since $\sss'^{\perp}$ is closed under extensions, it suffices to note that $\zzz = \sss^{\perp} \subset {\sss'}^{\perp}$ is the induced deformation of $\sss^{\perp}$. For the other inclusion, consider $X \in \sss'^{\perp}$.
By \cite[Prop. 4.2.1]{lowen5} and symmetry, $\zzz' = \{ C \in \ddd\, |\, \Hom_S(R,C) \in \zzz\}$, so we are to show that $\Hom_S(R,X) \in \zzz = \sss^{\perp}$. This follows from the way $\sss^{\perp} \subset {\sss'}^{\perp}$ is induced.
(2) Suppose $\sss$ and $\zzz$ are compatible. By Proposition \ref{compprop} and symmetry, to show that $\sss'$ and $\zzz'$ are compatible it suffices to show that $q_{\zzz'}(\sss') \subset \sss'$ for $q_{\zzz} = i_{\zzz'}a_{\zzz'}$ associated to $\zzz'$. So consider $X \in \sss'$. It suffices that $\Hom_S(R, q_{\zzz'}(X)) \in \sss$. We claim that $\Hom_S(R, q_{\zzz'}(X)) = q_{\zzz}(\Hom_S(R,X))$, which is in $\sss$ by assumption. We compute that $j\Hom_S(R, q_{\zzz'}(X)) =q_{\zzz'}(jHom_S(R,X)) = jq_{\zzz}(\Hom_S(R,X))$ where we have used that $q_{\zzz'}$ is left exact.
Now it remains to show that $(\sss' \ast \zzz')\cap \ccc = \sss \ast \zzz$. It is equivalent to show that $\sss' \ast \zzz'$ is the smallest Serre subcategory containing $S \ast T$. Since any Serre subcategory containing $S$ and $T$ also contains $\sss'$ and $\ttt'$ and hence $\sss' \ast \ttt'$, the proof is complete.
\end{proof}

\begin{proposition}
For a cover $U_i \lra U$, the canonical $\Des(U_i,\ccc) \lra \Des(U_i,\ddd)$ is a deformation fitting into a diagram
$$\xymatrix{{\ddd(U)} \ar[r] &{\Des(U_i,\ddd)} \\ {\ccc(U)}\ar[u] \ar[r] &{\Des(U_i,\ccc)}\ar[u]}$$ 
\end{proposition}

\begin{proof}
This follows from the fact that colimits (and hence also $R$-objects) in $\Des(U_i,\ddd)$ are pointwise.
\end{proof}

\begin{proposition}
Let $\fff$ be a fibered category of localizations and let $U_i \lra U$ be a cover.
\begin{enumerate}
 \item The category $\Des(U_i,\fff)$ is a cocomplete abelian category with exact filtered colimits.
 \item The functor $\fff(U) \lra \Des(U_i,\fff)$ is exact and colimit preserving.
 \item The functor $\fff(U) \lra \Des(U_i,\fff)$ has a right adjoint $l:  \Des(U_i,\fff) \lra \fff(U)$.
 \end{enumerate}
\end{proposition}

\begin{proof}
Since the restriction functors are exact left adjoints, it is easily seen that all finite limit and arbitrary colimit constructions can be carried out ``pointwise'' in $\Des(U_i,\fff)$, hence (1) and (2) follow. For the definition of $l$, consider a descent datum $F_V \in \fff(V)$ for $V \subset U_i$ with isomorphisms $F_V|_W \cong F_W$ for $W \subset V$. This determines a diagram in $\fff(U)$ consisting of the maps $j_{U,V\ast}F_V \lra j_{U,W\ast}F_W$. We define $l(F_V)$ to be the limit of this diagram in $\fff(U)$. Consequently, for $F \in \fff(U)$, we have $\Hom_{\fff(U)}(F, l(F_V)) = \lim \Hom_{\fff(U)}(F,j_{U,V\ast}F_V) = \lim \Hom_{\fff(U_V)}(F|_{V}, F_V) = \Hom_{\Des(U_i,\fff)}(F|_{V},F_V)$.
\end{proof}

\begin{proposition}\label{stacklift}
Suppose $\ccc$ is a fibered category of localizations and $\ddd$ is a flat localization deformation of $\ccc$. Suppose $U_i \lra U$ is a covering for which  the category $\Des(U_i,\ddd)$ is flat. If $\ccc(U) \lra \Des(U_i,\ccc)$ is an equivalence, then so is $\ddd(U) \lra \Des(U_i,\ddd)$. In particular, if $\ccc$ is a stack and all the categories $\Des(U_i,\ddd)$ are flat, then $\ddd$ is a stack, hence a localization stack deformation of $\ccc$.
\end{proposition}

\begin{proof}
The exact functor $\ddd(U) \lra \Des(U_i,\ddd)$ preserves coflat objects (\cite{lowenvandenbergh1}), and its right adjoint preserves injectives.
Consequently, we can apply \cite[Theorem 7.3]{lowenvandenbergh1} twice to obtain that both functors are fully faithful, hence they constitute an equivalence.
\end{proof}

We will now point out two situations in which the previous theorem applies.

\begin{proposition}\label{desflat}
Consider an open covering $U_i \lra U$  of $U$ which is closed under intersections. Suppose $\ccc$ is a flat fibered category of localizations on the $U_i$. In any of the following two cases, The category $\Des(U_i, \ccc)$ is flat:
\begin{enumerate}
\item If $\ccc$ is complemented, then the restriction functor $\Des(U_i, \ccc) \lra \ccc(U_i)$ has an exact left adjoint and the category $\Des(U_i, \ccc)$ is flat.
\item If the collection of $U_i$ is finite and the functors $j_{U_k,U_i,\ast}$ are exact, then the restriction functor $\Des(U_i, \ccc) \lra \ccc(U_i)$ has an exact right adjoint and the category $\Des(U_i, \ccc)$ is flat.
\end{enumerate}
\end{proposition}

\begin{proof}
(1) The restriction functor $j^{\ast}_{U_i}: \Des(U_i, \ccc) \lra \ccc(U_i)$ has an exact left adjoint $j_{U_i!}$ defined by $$j^{\ast}_{U_k}j_{U_i!}C_i = j_{U_k, U_i \cap U_k!} j^{\ast}_{U_i,U_i \cap U_k}C_i$$ Note that this defines a descent datum by compatibility (see Proposition \ref{complift}). Consider $C \in \Des(U_i, \ccc)$. We have an epimorphism $\coprod_i j_{U_i!}j_{U_i}^{\ast}C \lra C$. For every $i$ we pick an epimorphic effacement $P_i \lra j_{U_i}^{\ast}C$ for $\Tor_1(X,-)$ in $\ccc(U_i)$. Since $j_{U_i!}$ is exact and $\coprod$ is exact in $\Des(U_i, \ccc)$, $\coprod_i j_{U_i!}P_i \lra \coprod_i j_{U_i!}j_{U_i}^{\ast}C$ is an epimorphic effacement for $\Tor_1(X,-)$ in $\Des(U_i, \ccc)$.
(2)  The restriction functor $j^{\ast}_{U_i}: \Des(U_i, \ccc) \lra \ccc(U_i)$ has an exact right adjoint $j_{U_i\ast}$ defined by $$j^{\ast}_{U_k}j_{U_i\ast}C_i = j_{U_k, U_i \cap U_k\ast} j^{\ast}_{U_i,U_i \cap U_k}C_i$$ Note that this defines a descent datum by compatibility (see Proposition \ref{complift}). Consider $C \in \Des(U_i, \ccc)$. We have a monomorphism $C \lra \oplus_i j_{U_i\ast}j_{U_i}^{\ast}C$. For every $i$ we pick a monomorphic effacement $j_{U_i}^{\ast}C \lra E_i$ for $\Ext^1(X,-)$ in $\ccc(U_i)$. Since $j_{U_i\ast}$ is exact and $\oplus$ is exact in $\Des(U_i, \ccc)$, $\oplus_i j_{U_i\ast}j_{U_i}^{\ast}C \lra \oplus_i  j_{U_i\ast} E_i$ is a monomorphic effacement for $\Ext^1(X,-)$ in $\Des(U_i, \ccc)$.
\end{proof}

\subsection{Lifting (G), (F) and (FF) under deformation} \label{liftGFFF}

This section contains a technical result on how the conditions (G), (F) and (FF) (see \S \ref{parGFFF}) can be lifted under deformation. This result (Theorem \ref{liftGFFFthm}) will be used in the main Theorem \ref{theorem} of the main section \S \ref{abdefstack}.

Consider a diagram
$$\xymatrix{{\UUU} \ar[r]^u \ar[d]_f & {\ccc} \ar[d]^{S \otimes_R -} \\
{\VVV} \ar[r]_v & {\ddd}}$$
in which
\begin{itemize}
\item $u$ is an $R$-linear functor from a small $R$-linear category to an $R$-linear Grothendieck category.
\item $v$ is an $S$-linear functor from a small $S$-linear category to an $S$-linear Grothendieck category.
\item $u(\UUU)$ consists of flat objects.
\item $f$ is a \emph{not necessarily flat} deformation of linear categories.
\item $S \otimes_R -$ is left adjoint to an abelian deformation. 
\end{itemize}

\begin{theorem}\label{liftGFFFthm}
If $v$ satisfies (G), (F) and (FF), then so does $u$.
\end{theorem}

\begin{proof}[Proof of (F) and (FF)]
Consider the following diagram of presheaves on $\UUU$:
$$\xymatrix{{\Tor_1^R(S,\UUU(-,U))} \ar[r] \ar[d]  & {0} \ar[d]\\
{I \otimes_R \UUU(-,U)} \ar[r] \ar[d] & {\ccc(u(-), I \otimes_R u(U))} \ar[d]\\
{\UUU(-,U)} \ar[r] \ar[d] & {\ccc(u(-), u(U))} \ar[d] \\
{S \otimes_R \UUU(-,U)} \ar[r] \ar[d] & {\ccc(u(-), S \otimes_R u(U))} \ar[d] \\
{0} \ar[r] & {\Ext^1_{\ccc}(u(-), I \otimes_R u(U))}}$$
We are to prove that the middle arrow satisfies (F) and (FF). Since the lower one satisfies (FF) and the upper one satisfies (F), it suffices to show that for $X \in \mmod(S)$
$$X \otimes_S \VVV(f(-), f(U)) \lra \ddd(vf(-), X \otimes_S  vf(U)$$
satisfies (F) and (FF). Since both $f$ and $v$ satisfies (G) and (F), it suffices that
$$X \otimes_S \VVV(-,f(U)) \lra X \otimes_S \ddd(v(-),vf(U))$$ and  $$X \otimes_S \ddd(-,vf(U)) \lra \ddd(-, X \otimes_R vf(U))$$
satisfy (F) and (FF). The first map obviously does since $v$ satisfies (F) and (FF). The proof for the second map is similar to the proof of Proposition \ref{Xlemm}, but this time we use the fact that $\Ext^1_{\ddd}(-, D)$ is weakly effaceable.
\end{proof}

\begin{proof}[Proof of (G)]
Consider $C \in \ccc$, and write $C$ as an extension $E: 0 \lra IC \lra C \lra S \otimes_R C \lra 0$ of $\ddd$-objects. We will construct an epimorphism $\rho: \coprod_i u(U_i) \lra S \otimes_R C$ along which the pullback of $E$ splits. Consequently, $\rho$ will lift to an epimorphism $\coprod_i u(U_i) \lra C$.
First, take an epimorphism $\coprod_i u(U_i) \lra \coprod_i v(V_i) \lra S \otimes_R C$ where $f(U_i) = V_i$.
The pullback of $E$ is an extension $E': 0 \lra IC \lra C' \lra \coprod_i u(U_i)$. Since $\coprod_iu(U_i)$ is flat, $E'$ is the pullback of an extension $E'': 0 \lra IC \lra D \lra \coprod_i v(V_i)$ in $\ddd$. If we can split this extension by pulling back along a map $\coprod_j v(V_j) \lra \coprod_i v(V_i)$ which lifts to $\coprod_j u(U_j) \lra \coprod_i u(U_i)$, we are finished. For every finite sub coproduct $\coprod_K$ of $\coprod_i v(V_i)$, we consider the pullback $D_K$ of $D$. Since $v$ satisfies (G) and (F), we can generate $D_K$ with maps $v(V_{K,j}) \lra D_K$ such that all compositions $v(V_{K,j}) \lra v(V_i)$ are in the image of $v$. The induced map $\coprod v(V_{K,j}) \lra \coprod v(V_i)$ is epimorphic, splits $E''$ and lifts to $u(\UUU)$.
\end{proof}

\subsection{Relation between linear and abelian deformations}\label{abdefstack}

In this section, we will prove the following theorem:

\begin{theorem}\label{theorem}
Let $\aaa$ be a locally flat linear algebroid prestack on $X$. There is an equivalence between
\begin{enumerate}
\item locally flat weak linear deformations of $\aaa$
\item flat abelian deformations of the abelian category $\Mod(\aaa)$
\item flat localization deformations of the fibered category $\rMod(\aaa)$
\item flat localization stack deformations of the stack $\rMod(\aaa)$
\end{enumerate}
If $\uuu \subset \mathrm{Open}(X)$ is a \emph{basis} of $X$, there is a further equivalence with
\begin{enumerate}
\item[(5)] locally flat weak linear deformations of $\aaa|_{\uuu}$
\end{enumerate}
\end{theorem}

\begin{proof}
By Proposition \ref{largest}, there is an equivalence between (2) and (3). If $\ddd$ is a flat localization deformation of $\rMod(\aaa)$, then by Proposition \ref{complemented}, $\ddd$ is complemented, so by Propositions \ref{desflat}, \ref{stacklift}, $\ddd$ is a stack deformation, hence (3) and (4) coincide. It now suffices to show the equivalence of (2) and (5). Since $\uuu$ is a basis of $X$, we have an equivalence of categories $\Mod(\aaa|_{\uuu}) \cong \Mod(\aaa)$. Suppose $\bbb \lra \aaa|_{\uuu}$ is a locally flat weak linear deformation. 
By Proposition \ref{linnaarab}, there is an induced flat abelian deformation $\Mod(\aaa) \cong \Mod(\aaa|_{\uuu}) \lra \Mod(\bbb)$. Conversely, if $\Mod(\aaa) \lra \ddd(X)$ is a flat abelian deformation, we look at the induced localization stack deformation $\rMod(\aaa) \lra \ddd'$ on $X$. By Proposition \ref{abnaarlin} (1), there is an induced weak linear deformation $\aaa \lra \bbb$ which restricts to $\aaa|_{\uuu} \lra \bbb|_{\uuu}$. By Proposition \ref{algebroidclue} making a loop in (5) yields an equivalent deformation. By Proposition \ref{abnaarlin} (2), making a loop in (2) yields an equivalent deformation. This finishes the proof.
\end{proof}

\begin{proposition}\label{linnaarab}
Let $\uuu$ be a basis of $X$ and $\aaa$ a locally flat prestack on $\uuu$.
Suppose we are in one of the following cases:
\begin{enumerate}
\item $\bbb \lra \aaa$ is a flat linear deformation of fibered categories
\item $\bbb \lra \aaa$ is a locally flat weak linear deformation of prestacks
\end{enumerate}
There is an induced flat deformation of abelian categories $\Mod(\aaa) \lra \Mod(\bbb)$ and an induced flat deformation of  stacks $\rMod(\aaa) \lra \rMod(\bbb)$. The latter can be obtained from the former by taking induced deformations of the localizations $\Mod(\aaa_U)$ of $\Mod(\aaa)$ for $U \in \uuu$ . Moreover, equivalent deformations of prestacks yield equivalent abelian deformations.
\end{proposition}

\begin{proof}
The statement concerning flatness follows from Proposition \ref{platjes}.
If $\bbb \lra \aaa$ induces $\Mod(\aaa) \lra \Mod(\bbb)$, then $\bbb_U \lra \aaa_U$ induces $\Mod(\aaa_U) \lra \Mod(\bbb_U)$ in the same way, yielding a deformation of stacks which is indeed the induced stack of localizations.
In case (2), $R \otimes_S \bbb \lra \aaa$ consists of morphisms satisfying (G), (F) and (FF), so there is an induced equivalence ${\Mod}(\aaa) \cong {\Mod}(R \otimes_S \bbb)$ by Propositions \ref{assocmod}, \ref{GFFFstackequiv}.
So it remains to prove case (1).  Consider ${\Mod}(R \otimes_S \bbb) \lra {\Mod}(\bbb)$. Let $\BBB$ be the additive category associated to $\bbb$.
The additive category associated to $R \otimes_S \bbb$ is $R \otimes_S \BBB$.
By Proposition \ref{modsh}, we can look at $\Sh(R \otimes_S \BBB, \ttt_{R \otimes_S \bbb}) \lra \Sh(\BBB, \ttt_{\bbb})$ instead. By Proposition \ref{topos} and \cite{lowenvandenbergh1} \S 7, this map is a deformation if $\ttt_{\bbb} = (\ttt_{R \otimes_S \bbb})_{R \otimes_S -}$. Let $b_i: B_i \lra B|_{U_i}$ be a potential covering of $B \in \bbb(U)$ for $\ttt_{\bbb}$. Its image in $R \otimes_S \BBB$ is $R \otimes_S b_i: R \otimes_S B_i \lra (R \otimes_S B)|_{U_i}$. Since $R \otimes_S -: \bbb(U_i) \lra R \otimes_S \bbb(U_i)$ reflects isomorphisms, $b_i$ is covering for $\ttt_{\bbb}$ if and only if $R \otimes_S b_i$ is covering for $\ttt_{R \otimes_S \bbb}$, which finishes the proof.
\end{proof}

We have used the following relation between topologies on $\BBB$ and $R \otimes_S \BBB$ for a linear category $\BBB$:

\begin{proposition}\label{topos}
Let $\BBB$ be an $S$-linear category. Consider $R \otimes_S -: \BBB \lra R \otimes_S \BBB$. 
\begin{enumerate}
\item If $\ttt$ is a topology on $\BBB$, then $R \otimes_S \ttt$ is a topology on $R \otimes_S \BBB$.
\item If $\ttt$ is a topology on $R \otimes_S \BBB$, then $\ttt_{R \otimes_S -}$ is a topology on $\BBB$.
\end{enumerate}
Taking generated subfunctors yields a one one correspondence between Grothendieck topologies on $\BBB$ and $R \otimes_S \BBB$, compatible with the one one correpondence between localizations of $\Mod(\BBB)$ and $\Mod(R \otimes_S \BBB)$ of \cite{lowenvandenbergh1}, \S 7.
\end{proposition}

\begin{proof}
This is easily deduced from the correspondence between localizing subcategories in $\Mod(\BBB)$ and $\Mod(R \otimes_S \BBB)$ of \cite{lowenvandenbergh1}, \S 7.
\end{proof}

We will now give a converse to Proposition \ref{linnaarab}.

\begin{proposition}\label{abnaarlin}
Let $\aaa$ be a locally flat $R$-linear prestack on $X$. Suppose that $\rMod(\aaa) \lra \ddd$ is a flat localization stack $S$-deformation. Let $\bbb = \bar{\aaa}$ be the full sub-prestack of $\ddd$ spanned by the flat objects $D \in \ddd(U)$ with $R \otimes_S D \in \aaa(U)$. We have morphisms
$$\xymatrix{{\bbb} \ar[r] \ar[d] & {\ddd} \ar[d] \\ {\aaa} \ar[r] & {\rMod(\aaa)}}$$
The following is true:
\begin{enumerate}
\item $\bbb \lra \aaa$ is a locally flat weak linear deformation of prestacks
\item $\bbb \lra \ddd$ yields an equivalence $\rMod(\bbb) \cong \ddd$
\end{enumerate}
\end{proposition}

\begin{proof}

(1) follows from Proposition \ref{defalg}.
To prove (2), we use Theorem \ref{karmodthm}. First note that by Proposition \ref{complemented}, $\ddd$ is a complemented stack of localizations, so the theorem applies. To prove part (a), let $\BBB$ be the additive category of $\bbb$. By (1) we have a diagram
$$\xymatrix{{{\BBB}_U} \ar[r] \ar[d] & {\ddd(U)} \ar[d] \\ {R \otimes_S {\BBB}_U} \ar[r] & {\Mod(\aaa|_U)}}$$
in which the lower arrow satisfies (G), (F) and (FF). It follows from \S \ref{liftGFFF}, Theorem \ref{liftGFFFthm} that the upper arrow also satisfies (G), (F) and (FF). 
Next we prove part (b). On $\AAA_U$, by assumption, the covering systems $\ttt_{\mathrm{epi},U}$ and $\ttt_U$ are equivalent. We have to show that the same holds on $\BBB_U$. Consider $\varphi: \BBB_U \lra \AAA_U$ and consider a $\ttt_{\mathrm{epi},U}$-covering $b_i:B_i \lra B$ in $\BBB_U$. The images $\varphi(b_i):\varphi(B_i) \lra \varphi(B)$ constitute a $\ttt_{\mathrm{epi},U}$-covering of $\varphi(B)$ in $\AAA_U$. Consequently, there is a covering $W_k \lra V$ such that for the canonical $b_k: B|_{W_k} \lra B$,  the morphisms $\varphi(b_k)$ are in the generated subfunctor, i.e. there are finitely many $c_i: \varphi(B|_{W_k}) \lra \varphi(B_i)$ with $\sum_i \varphi(b_i)c_i = \varphi(b_k)$. The morphisms $c_i$ are not in the image of $\varphi$, but by Proposition \ref{defalg}, for every $k$ there is a  covering $W_{kl} \lra W_k$ such that for the canonical $b_{kl}: B|_{W_{kl}} \lra B|_{W_k}$ we have $c_i\varphi(b_{kl}) = \varphi(b_{ikl})$. Consequently, $\varphi(\sum_i b_i b_{ikl}) = \varphi(b_k b_{kl})$, so the morphism $\sum_i b_i b_{ikl}$, which is in the subfunctor generated by the covering $b_i$, is defined by some $\bbb(W_{kl})$-morphism $b: B|_{W_{kl}} \lra B|_{W_{kl}}$, which is mapped by $\bbb(W_{kl}) \lra \aaa(W_{kl})$ onto an isomorphism. Since $\bbb(W_{kl})$ is a full subcategory of flat objects in $\ddd(W_{kl})$, it follows that $b$ is itself an isomorphism, which finishes the proof.
\end{proof}

\begin{proposition}\label{algebroidclue}
Let $\bbb \lra \aaa$ be a locally flat weak linear deformation of algebroid prestacks on a basis $\uuu$ of $X$. 
We have morphisms
$$\xymatrix{{\bbb} \ar[r] \ar[d] & {\bar{\aaa}}\ar[r] \ar[d] & {\rMod(\bbb)} \ar[d] \\ {\aaa} \ar[r] & {\aaa} \ar[r] & {\rMod(\aaa)}}$$
in which $\bar{\aaa} \lra \aaa$ is the deformation constructed in Proposition \ref{abnaarlin} (restricted to $\uuu$). The morphism $\bbb \lra \bar{\aaa}$ is fully faithful and satisfies (G).
\end{proposition}

\begin{proof}
This follows from Proposition \ref{algG}.
\end{proof}

The remainder of this section contains auxiliary results for the proof of Proposition \ref{abnaarlin}. We use the notations of Proposition \ref{abnaarlin}.
The following key-result, which proves Proposition \ref{abnaarlin}(1), makes use of the obstruction theory for lifting objects and maps in a deformation of abelian categories, as developed in \cite{lowen2}.

\begin{proposition}\label{defalg}
$\bbb \lra \aaa$ satisfies (G) and (F) and $R \otimes_S \bbb \lra \aaa$ satisfies (FF).
\end{proposition}

\begin{proof}
For $A_U \in \aaa(U)$, the obstruction against lifting $A_U$ is $o \in \Ext^2(A_U, I \otimes_R A_U)$. By local projectivity of $A_U$, we can take a cover $U_i \lra U$ such that $o|_{U_i} = 0$, hence the objects $A_U|_{U_i}$ lift to $\ddd(U_i)$. For a map $R \otimes_S B_U \lra R \otimes_S B_U'$, the obstruction against lifting is in $\Ext^1(R \otimes_S B_U, I \otimes_S B_U')$, so the same argument applies.  
Now consider $b: B \lra B'$ with $R \otimes_S b = 0$. So $b \in \ddd(U)(B, I \otimes_S B')$. By Propositions \ref{Extnul} and \ref{Xlemm}, there is a cover $U_i \lra U$ for which $b|_{U_i}$ is in the image of $I \otimes_S \ddd(U_i)(B|_{U_i}, B'|_{U_i})$. Consequently, $b|_{U_i}$ is zero in $R \otimes_S \bbb(U_i)(B|_{U_i}, B'|_{U_i})$, as desired.
 \end{proof}
 
 \begin{proposition}\label{Xlemm}
 Suppose $B, B'$ are flat in $\ddd(U)$ and $B$ is locally projective. For $X \in \mmod(R)$, the morphism of \emph{presheaves} on $\uuu/U$
 $$X \otimes_R \Hom_{\ddd}(B,B') \lra \Hom_{\ddd}(B, X \otimes_R B')$$
 satisfies (F) and (FF).
\end{proposition}

\begin{proof}
Write $0 \lra K \lra R^n \lra X \lra 0$ and consider the following diagram:
$$\xymatrix{{K \otimes_R \Hom_{\ddd}(B,B')} \ar[r]_{\alpha_1} \ar[d] & {\Hom_{\ddd}(B, K \otimes_R B')} \ar[d]\\
{R^n \otimes_R \Hom_{\ddd}(B,B')} \ar[r]_{\alpha_2}^{\cong} \ar[d] & {\Hom_{\ddd}(B, R^n \otimes_R B')} \ar[d]\\
{X \otimes_R \Hom_{\ddd}(B,B')} \ar[r]_{\alpha_3} \ar[d] & {\Hom_{\ddd}(B, X \otimes_R B')} \ar[d]\\
{0} \ar[r]_{\alpha_4} \ar[d] & {\Ext^1_{\ddd}(B, K \otimes_R B')} \ar[d] \\
{0} \ar[r]_{\alpha_5}  & {\Ext^1_{\ddd}(B, R^n \otimes_R B')}}$$
Since $\alpha_5$ satisfies (FF) and both $\alpha_2$ and $\alpha_3$ satisfy (F), 
$\alpha_3$ satisfies (F). Consequently, $\alpha_1$ too satisfies (F). Hence, since both $\alpha_2$ and $\alpha_4$ satisfy (FF), the same holds for $\alpha_3$.
\end{proof}

\begin{proposition}\label{Extnul}
Let $B \in \ddd(U)$ be a flat object. If $R \otimes_S B$ is locally projective (resp. locally finitely presented) in $\rMod(\aaa)$, the same holds for $B$ in $\ddd$.
\end{proposition}

\begin{proof}
For $M \in \Mod(\aaa|_U)$, $\Ext^i_{\ddd}(B,M) = \Ext^i_{\rMod(\aaa)}(R \otimes_S B,M)$. For general $M$, it suffices to write $M$ as an extension of objects in $\Mod(\aaa|_U)$.
\end{proof}

\subsection{The case of an acyclic basis}\label{paracyclic}
In this section we will briefly discuss the relationship with some results in \cite{lowenvandenbergh1}. In particular we will reprove \cite[Theorem 8.18]{lowenvandenbergh1}.
We consider an $R$-linear algebroid prestack $\aaa$ on $X$ and we suppose that $X$ has an \emph{acyclic basis} $\uuu$ for $\aaa$, i.e. for every $U \in \uuu$, $A, A' \in \aaa(U)$, $X \in \mmod(R)$, $i = 1,2$ we have
$$\Ext^i_{\Mod(\aaa|_U)}(A, X \otimes_R A') = 0$$
By Theorem \ref{theorem}, deforming the abelian category $\Mod(\aaa)$ is equivalent to weakly deforming the prestack $\aaa|_{\uuu}$. In this section we show that it is also equivalent to deforming $\aaa|_{\uuu}$ \emph{as a fibered category}. The following Theorem essentially generalizes \cite[Theorem 8.18]{lowenvandenbergh1}.
\begin{theorem}\label{acyc}
Let $\aaa$ be a locally flat $R$-linear algebroid prestack on $X$ and let $\uuu$ be an acylic basis for $\aaa$. Then $\aaa|_{\uuu}$ is flat as a fibered category and there is an equivalence between:
\begin{enumerate}
\item flat deformations of the fibered category $\aaa|_{\uuu}$
\item flat deformations of the linear category $\AAA_{\uuu}$ associated to $\aaa|_{\uuu}$
\item flat abelian deformations of the abelian category $\Mod(\aaa)$
\end{enumerate}
\end{theorem} 

\begin{proof}
The equivalence of (1) and (2) follows from Propositions \ref{gradlin} and \ref{gradfib}.
The remainder of the proof is a modification of the proof of Theorem \ref{theorem}.
Let $\bbb \lra \aaa|_{\uuu}$ be a flat linear deformation of fibered categories. By Proposition \ref{linnaarab}, this yields a flat abelian deformation $\Mod(\aaa) \cong \Mod(\aaa_{\uuu}) \lra \Mod(\bbb)$. Conversely for a flat abelian deformation $\Mod(\aaa) \lra \ddd(X)$ we look at the induced localization stack deformation $\rMod(\aaa) \lra \ddd$. This yields a weak deformation $\bar{\aaa} \lra \aaa$ as in Proposition \ref{abnaarlin}, only this time, by acyclicity, every $A \in \aaa(U)$ with $U \in \uuu$ has (up to isomorphism) a unique flat lift to $\ddd(U)$, and every $\bar{\aaa}(U) \lra \aaa(U)$ is actually a deformation of linear categories. It remains to say why the linear deformation $\bbb \lra \aaa|_{\uuu}$ is equivalent to $\bar{\aaa}|_{\uuu} \lra \aaa|_{\uuu}$ obtained from $\rMod(\aaa|_{\uuu}) \lra \rMod(\bbb)$. Certainly $\bbb \lra \rMod(\bbb)$ factors over $\bbb \lra \bar{\aaa}|_{\uuu}$ and the functors $\bbb(U) \lra \bar{\aaa}(U)$ are essentially surjective. That they are also fully faithful easily follows from the 5-lemma since the functors $\aaa(U) \lra \Mod(\aaa|_U)$ are.
\end{proof}

\begin{corollary} Let $\ooo$ be a locally flat (twisted) sheaf of $R$-algebras on $X$ and let $\uuu$ be an acyclic basis for $\ooo$. Then $\ooo|_{\uuu}$ is flat as a (twisted) presheaf and there is an equivalence between:
\begin{enumerate}
\item flat deformations of $\ooo|_{\uuu}$ as a twisted presheaf
\item flat deformations of the linear category $\OOO_{\uuu}$ associated to $\ooo|_{\uuu}$
\item flat deformations of the abelian category $\Mod(\ooo)$
\end{enumerate}
\end{corollary}

\begin{proof}
This follows from Theorem \ref{acyc} and Proposition \ref{twistfibgrad}.
\end{proof}

\begin{remark}\label{remacyc}
By \cite[Prop. 6.13]{lowenvandenbergh1}, if $\bbb$ is a locally flat weak linear deformation of $\aaa$ and $U$ is an acyclic open subset for $\aaa$, then $U$ is an acyclic open subset for $\bbb$. Hence for algebroid prestacks on $X$, an acyclic basis ``lifts'' under flat weak linear deformation.
\end{remark}

\subsection{Quasi-coherent sheaves}\label{parafcover}

If $(X,\ooo)$ is a quasi-compact, separated scheme, we know that the abelian categories $\Mod(\ooo)$ of sheaves of modules on $X$, and $\Qch(\ooo)$ of quasi-coherent sheaves of modules on $X$ have the same Hochschild cohomology \cite{lowenvandenbergh2}. In this section we prove the deformation analogon of this result.

Let $\aaa$ be an algebroid prestack on $X$. The restriction morphisms $\aaa(U) \lra \aaa(V)$ can be uniquely extended to colimit-preserving functors $\Mod(\aaa(U)) \lra \Mod(\aaa(V))$ defining a fibered category of abelian categories which we will denote $\aaa dd(\aaa)$ to distinguish in notation from $\rMod(\aaa)$.
If $\uuu$ is a covering of $X$, we put $\Qch(\uuu, \aaa) = \Des(\uuu, \aaa dd(\aaa))$. We consider the fibered category $\qqq ch(\uuu, \aaa)$ on $\uuu \cup \{X\}$ with $\qqq ch(\uuu, \aaa)(U) = \Mod(\aaa(U))$ and $\qqq ch(\uuu, \aaa)(X) = \Qch(\uuu, \aaa)$.

\begin{proposition}\label{descong}
Let $\aaa$ and $\uuu$ be as above and suppose $\vvv \subset \uuu$ is a finite covering of $X$ which is closed under intersections, and such that every $U \in \uuu$ is contained in some $V \in \vvv$. Then we have an equivalence of categories
$\Qch(\uuu, \aaa) \cong \Qch(\vvv, \aaa)$.\qed
\end{proposition}

\begin{theorem}\label{quasi}
Let $\aaa$ be a locally flat $R$-linear algebroid prestack on $X$ with an acyclic basis $\uuu$. Suppose $\vvv \subset \uuu$ is a finite covering of $X$ which is closed under intersections and such that every $U \in \uuu$ is contained in some $V \in \vvv$.
Suppose that $\qqq ch(\uuu, \aaa)$ is a fibered category of localizations.
The fibered category $\aaa|_{\uuu}$ is flat and the abelian category $\Qch(\uuu, \aaa) \cong \Qch(\vvv, \aaa)$ is flat.
There is an equivalence between: 
\begin{enumerate}
\item flat abelian deformations of the abelian category $\Mod(\aaa)$
\item flat linear deformations of the fibered category $\aaa|_{\vvv}$
\item flat abelian deformations of the abelian category $\Qch(\uuu, \aaa)$
\end{enumerate}
\end{theorem}

\begin{proof}
By Theorem \ref{acyc}, (1) can be replaced by 
\begin{enumerate}
\item[(1')] flat linear deformations of the fibered category $\aaa|_{\uuu}$.
\end{enumerate}
Flatness of $\aaa|_{\uuu}$ is stated in Theorem \ref{acyc}. By Proposition \ref{descong}, $\Qch(\uuu, \aaa) \cong \Qch(\vvv, \aaa)$.
Consider the fibered category $\aaa dd(\aaa)|_{\vvv}$. For $V' \subset V$ in $\vvv$, the right adjoint of $\Mod(\aaa(V)) \lra \Mod(\aaa(V'))$ is the forgetful functor composition with $\aaa(V) \lra \aaa(V')$, which is obviously exact. Hence by Proposition \ref{desflat}(2), the category $\Qch(\vvv, \aaa) = \Des(\vvv, \aaa dd(\aaa))$ is flat.

We can go from (1') to (2) simply by restriction. Suppose we have a flat linear deformation $\bbb \lra \aaa|_{\vvv}$. The corresponding abelian deformations $\Mod(\aaa(V)) \lra \Mod(\bbb(V))$ obviously define a flat abelian deformation $\aaa dd(\aaa|_{\vvv}) \lra \aaa dd(\bbb)$, which is a localization deformation by \cite[Theorem 7.3]{lowenvandenbergh1}. Again by Proposition \ref{desflat}(2), the category $\Des(\vvv, \aaa dd(\bbb))$ is a flat abelian deformation of $\Qch(\vvv, \aaa)$, and we arrive at (3). Also, the restrictions $\Des(\vvv, \aaa dd(\bbb)) \lra \Mod(\bbb(V))$ have exact adjoints, so by \cite[Theorem 7.3]{lowenvandenbergh1} they constitute localizations.

Now suppose $\Qch(\uuu, \aaa) \lra \qqq(X)$ is a flat abelian deformation. The fibered category of localizations $\qqq ch(\uuu, \aaa)$ on $X$ induces a fibered category of localizations $\qqq$ deforming $\qqq ch(\uuu, \aaa)$. For $U \in \uuu$, the projective generators $\aaa(U)$ of $\Mod(\aaa(U))$ can be uniquely lifted to $\qqq(U)$, yielding equivalences $\Mod(\bbb(U)) \lra \qqq(U)$ and functors $\bbb(U) \lra \bbb(U')$ turning $\bbb$ into a linear deformation of $\aaa|_{\uuu}$, and we arrive at (1'). 

To see that a loop in (3) yields equivalent deformations, it suffices to note that  by flatness of $\Des(\vvv, \aaa dd(\bbb))$, we get an equivalence $\qqq(X) \lra \Des(\vvv, \aaa dd(\bbb))$ by Proposition \ref{stacklift}. To see that a loop in (1') or (2) yields equivalent deformations, it suffices to note that for $U \in \uuu$, $U \subset V$, $V \in \vvv$, starting from (1') or (2) we get localizations $\Mod(\bbb(U)) \lra \Mod(\bbb(V)) \lra \Des(\vvv, \aaa dd(\bbb))$ which are then necessarily the induced deformations. 
\end{proof}

\begin{corollary}\label{quasicor}
Suppose $(X,\ooo)$ is a quasi-compact, separated scheme such that $\ooo$ is flat in $\Mod(\ooo)$. Then the abelian categories $\Mod(\ooo)$ and $\Qch(\ooo)$ are flat and there is an equivalence between:
\begin{enumerate}
\item flat abelian deformations of $\Mod(\ooo)$
\item flat abelian deformations of $\Qch(\ooo)$
\end{enumerate}
If $\vvv$ is a finite affine covering of $X$ closed under intersections, there is a further equivalence with
\begin{enumerate}
\item[(3)] flat deformations of $\ooo|_{\vvv}$ as a twisted presheaf
\item[(4)] flat deformations of the linear category $\OOO_{\vvv}$ associated to $\ooo|_{\vvv}$
\end{enumerate}
\end{corollary}

\begin{remark}
For an arbitrary algebroid prestack $\aaa$ on $X$, we can define the \emph{stack of quasi-coherent sheaves} over $\aaa$ to be the associated stack of $\aaa dd(\aaa)$, i.e. $\qqq ch(\aaa) = as(\aaa dd(\aaa))$. The abelian category of quasicoherent sheaves on an open subset $U$ is by definition the category $\qqq ch(\aaa)(U)$ and $\Qch(\aaa) = \qqq ch(\aaa)(X)$. If $(X, \ooo)$ is a scheme, then $\aaa dd(\ooo)|_{\uuu}$ is a stack on the basis $\uuu$ of affine opens. Consequently, $\Qch(\ooo) = \Qch(\uuu, \ooo)$. It is not (yet) clear to us in which generality there is an equivalence between the abelian deformations of $\Mod(\aaa)$ and of $\Qch(\aaa)$.
\end{remark}

\def\cprime{$'$}
\providecommand{\bysame}{\leavevmode\hbox to3em{\hrulefill}\thinspace}
\providecommand{\bysame}{\leavevmode\hbox to3em{\hrulefill}\thinspace}
\providecommand{\MR}{\relax\ifhmode\unskip\space\fi MR }
\providecommand{\MRhref}[2]{%
  \href{http://www.ams.org/mathscinet-getitem?mr=#1}{#2}
}
\providecommand{\href}[2]{#2}

\bibliographystyle{amsplain}
\end{document}